\documentclass[10pt, reqno]{amsart}
\usepackage{tikz}
\usetikzlibrary{matrix,arrows,decorations.pathmorphing}
\usepackage{amssymb}

\newtheorem{proposition}{Proposition}[section]
\newtheorem{lemma}[proposition]{Lemma}
\newtheorem{corollary}[proposition]{Corollary}
\newtheorem{theorem}[proposition]{Theorem}
\newtheorem{remark}[proposition]{Remark}

\theoremstyle{definition}

\newtheorem{example}[proposition]{Example}

\newcommand{\selabel}[1]{\label{se:#1}}

\newcommand{\eqlabel}[1]{\label{eq:#1}}
\newcommand{\equref}[1]{(\ref{eq:#1})}

\def\<{\leqslant}
\def\>{\geqslant}
\def\a{\alpha}
\def\b{\beta}

\def\g{\gamma}
\def\G{\Gamma}

\def\t{\triangle}

\def\l{\lambda}
\def\L{\Lambda}
\def\s{\sigma}

\def\ot{\otimes}

\date{}

\begin{document}
\title[The Projective Class Rings of Drinfeld doubles pointed rank one Hopf algebras]{The Projective Class Rings of  Drinfeld doubles of pointed rank one Hopf algebras}
\author{Hua Sun}
\address{School of Mathematical Science, Yangzhou University,
	Yangzhou 225002, China}
\email{huasun@yzu.edu.cn}
\author{Hui-Xiang Chen}
\address{School of Mathematical Science, Yangzhou University,
	Yangzhou 225002, China}
\email{hxchen@yzu.edu.cn}
\author{Libin Li}
\address{School of Mathematical Science, Yangzhou University,
	Yangzhou 225002, China}
\email{lbli@yzu.edu.cn}
\author{Yinhuo Zhang}
\address{Department of Mathematics $\&$  Statistics, University of Hasselt, Universitaire Campus, 3590 Diepenbeek, Belgium
	Yangzhou 225002, China}
\email{yinhuo.zhang@uhasselt.be}
\begin{abstract}
Let $\Bbbk$ be an algebraically closed field of characteristic $0$.
In this paper, we study the Grothendieck ring $G_0(D(H_\mathcal{D}))$ and the projective class ring $r_p(D(H_\mathcal{D}))$ of the Drinfeld double $D(H_{\mathcal{D}})$ of the rank one pointed Hopf algebra $H_{\mathcal{D}}$.
We analyze the tensor products of simple modules with simple modules, simple modules with indecomposable projective modules, and indecomposable projective modules with indecomposable projective modules, providing explicit decomposition rules in each case.
Finally, we compute both the Grothendieck ring $G_0(D(H_\mathcal{D}))$ and the projective class ring $r_p(D(H_\mathcal{D}))$, and present these two rings in terms of generators and defining relations.
\end{abstract}
\maketitle

\section{Introduction}\selabel{1}
The representation theory of finite-dimensional Hopf algebras occupies a central position in modern algebra, connecting deeply with quantum groups, tensor categories, and low-dimensional topology. Among the most fundamental algebraic invariants associated with such algebras $H$ are the Grothendieck ring $G_{0}(H)$  and the projective class ring $r_p(H)$, which encode the tensor product structures of the categories of simple and projective
$H$-modules, respectively.
These rings, often viewed as fusion rings of the corresponding tensor categories, provide valuable insight into how representations interact under tensor operations and how the underlying algebraic symmetries manifest at the categorical level.

The family of rank one pointed Hopf algebras, classified by Krop and Radford \cite{KropRad} over an algebraically closed field of characteristic zero, and by Scherotzke in a positive characteristic \cite{Sche},  provides a natural generalization of the classical Taft algebras and Radford algebras.  Their representation theory has been developed in detail by Wang, Li and Zhang \cite{WangLiZhang2014, WangLiZhang2016}.  Although the representation category $_{H_{\mathcal D}}\mathcal{M}$ of a  rank one pointed Hopf algebra $H_{\mathcal D}$ is not a braided tensor category,  both its representation ring and its Grothendieck ring  are commutative.

However,  the tensor category $_{H_{\mathcal D}}\mathcal{M}$ can not, in general,  be embedded as a full tensor subcategory of the Drinfeld center of $H_{\mathcal D}$.  This observation motivated a deeper  study of the representation theory of the Drinfeld double $D(H_{\mathcal D})$,   and in particular, of the relationship between the Grothendieck (and Green)  rings of $D(H_{\mathcal D})$ and  those of $H_{\mathcal D}$.

Drinfeld doubles of specific classes of Hopf algebras, such as Taft algebras, generalized Taft algebras, and Radford algebras, have been extensively investigated. Chen and his collaborators  computed the structures of the Green, Grothendieck, and projective class  rings of  Drinfeld doubles of Taft algebras  in~\cite{Ch5, ChHasSun, ChHasLinSun, SunHasLinCh}. Erdmann et al. investigated the representations and stable Green rings of Drinfeld doubles of generalized Taft algebras in \cite{EGST2006, EGST2019}, while Sun and Chen studied the representations of Drinfeld doubles of Radford Hopf algebras in \cite{SunChen}, revealing intricate fusion rules closely related to those of small quantum groups and modular tensor categories.  Furthermore,  Krop and Radford classified all simple and projective indecomposable modules of the Drinfeld double $D(H_{\mathcal{D}})$  when  $G(H_{\mathcal{D}})$, the group of group-like elements, is abelian \cite{KropRad}. More recently,  Sun, Chen and Zhang in \cite{SUN C Z} completed  the classification of all finite-dimensional indecomposable $D(H_{\mathcal{D}})$-modules.

Since determining the complete fusion rules for all indecomposable modules of a rank one pointed Hopf algebra is an intensive task, we divide our project into  two stages: (1) first, we compute the decomposition rules for the tensor products of simple and projective indecomposable modules, leading to explicit descriptions of the Grothendieck and projective class rings; and
(2) subsequently, we will determine the fusion rules involving non-projective, non-simple indecomposable modules, from which the full Green and stable Green rings can be derived.

In this paper, we carry out the first of these steps.
We determine the Grothendieck ring $G_0(D(H_{\mathcal D}))$ and the projective class ring $r_p(D(H_{\mathcal D}))$ for an arbitrary rank one pointed Hopf algebra $H_{\mathcal D}$
 whose group of group-like elements is abelian.
We explicitly describe the tensor product decompositions among simple and indecomposable projective $D(H_{\mathcal D})$-modules,  and we present  the rings $G_0(D(H_{\mathcal D}))$ and $r_p(D(H_{\mathcal D}))$ in terms of generators and relations.
Our results unify and generalize previously known computations for the Drinfeld doubles of Taft algebras and other special cases of rank-one pointed Hopf algebras, providing a coherent framework for understanding their representation-theoretic and categorical structures.

The paper is organized as follows.
Section~\ref{s2} recalls the definition of a group datum, the construction of pointed rank-one Hopf algebras $H_{\mathcal D}$, and their Drinfeld doubles $D(H_{\mathcal D})$.
In Section~\ref{s3}, we review the simple and projective indecomposable $D(H_{\mathcal D})$-modules.
Section~\ref{s400} analyzes the socles of tensor products of simple modules, while Section~\ref{s500} provides explicit tensor product decompositions among simple and projective modules.
Finally, Section~\ref{s700} describes the structures of the Grothendieck ring $G_0(D(H_{\mathcal D}))$ and the projective class ring $r_p(D(H_{\mathcal D}))$.

Throughout, let $\Bbbk$ be an algebraically closed field with char$\Bbbk=0$ and $\Bbbk^{\times}=\Bbbk\backslash\{0\}$. Unless
otherwise stated, all algebras and Hopf algebras are
defined over $\Bbbk$; all modules are finite dimensional and left modules;
dim and $\otimes$ denote ${\rm dim}_{\Bbbk}$ and $\otimes_{\Bbbk}$,
respectively. Let $\mathbb Z$ denote the set of all integers, ${\mathbb Z}_n:={\mathbb Z}/n{\mathbb Z}$
for an integer $n$, and let $\mathbb{N}$ denote all non-negative integers.
The references \cite{Ka, Mon, Sw} are basic references for the theory of Hopf algebras and quantum groups. The readers can refer \cite{ARS} for the representation theory of algebras.

\section{Pointed rank one Hopf algebras and their doubles}\label{s2}

In this section, we recall the pointed rank one Hopf algebras and  their Drinfeld doubles.

Let $0\neq q\in\Bbbk$. For any integer $n>0$, define $$(n)_q=1+q+\cdots+q^{n-1}.$$ Note that $(n)_q=n$ when $q=1$ and $(n)_q=\frac{q^n-1}{q-1}$ when $q\neq 1$. Define the $q$-factorial of $n$ by $(0)!_q=1$ and $$(n)!_q=(n)_q(n-1)_q\cdots(1)_q$$ for $n>0$, see \cite[p.74]{Ka}.

A quadruple $\mathcal{D}=(G, \chi, a, \a)$ is called a {\it group datum} if $G$ is a finite group, $\chi$ is a $\Bbbk$-valued character of $G$, $a$ is a central element of $G$ and $\a\in\Bbbk$ subject to $\chi^n=1$ or $\a(a^n-1)=0$, where $n$ is the order of $\chi(a)$. The group datum $\mathcal{D}$ is  of {\it nilpotent type} if $\a(a^n-1)=0$, and it is of {\it non-nilpotent type} if $\a(a^n-1)\neq 0$ and $\chi^n=1$. For any group datum $\mathcal{D}=(G, \chi, a, \a)$, Krop and Radford constructed an associated finite dimensional pointed rank one Hopf algebra $H_{\mathcal{D}}$ and classified such Hopf algebras. They also described the Drinfeld doubles $D(H_{\mathcal D})$ of $H_{\mathcal D}$, see \cite{KropRad}.

Let $\mathcal{D}=(G, \chi, a, \a)$ be a group datum. The Hopf algebra $H_{\mathcal{D}}$ is generated as an algebra by $G$ and an element $x$, subject to the group relations for $G$, and the additional relations: $$x^n=\a(a^n-1) \text{ and } xg=\chi(g) gx\ \text{ for all }g\in G.$$
The comultiplication $\t$ is given by $$\t(x)=x\ot a+1\ot x \text{ and } \t(g)=g\ot g, \ g\in G.$$
Then $H_{\mathcal{D}}$ has a $\Bbbk$-basis $\{gx^j|g\in G, 0\<j\<n-1\}$. If $\mathcal{D}$ is of {\it non-nilpotent type}, then $H_{\mathcal{D}}\cong H_{\mathcal{D'}}$  as Hopf algebras, where $\mathcal{D'}=(G, \chi, a, 1)$. Therefore, we always assume that $\a=1$ whenever $\mathcal{D}$ is of {\it non-nilpotent type}. For further details, one can refer to \cite{KropRad, Rad1975}.

Throughout the following, let $\mathcal{D}=(G, \chi, a, \a)$ be a group datum, and assume that $G$ is abelian. Set $\rho=\chi(a)$ and denote by $n$ the order of $\rho$. Let $\Gamma={\rm Hom}(G,\Bbbk^{\times})$ denote the group of $\Bbbk$-valued characters of $G$.

 Recall that the Drinfeld double $D(H_{\mathcal D})=H_{\mathcal D}^{*\rm cop}\bowtie H_{\mathcal D}$ is generated, as an algebra, by its two sub-Hopf algebra $H_{\mathcal D}^{*\rm cop}$ and $H_{\mathcal D}$, where  $H^*_{\mathcal D}$ is the dual Hopf algebra of $H_{\mathcal D}$. The algebra $H_{\mathcal D}^{*\rm cop}$ is generated  by $\xi$ and the group $\Gamma$, subject to the following relations:
$$\xi^n=0 \text{ and }  \xi\gamma=\gamma(a)\gamma\xi \text{ for all } \gamma\in\Gamma.$$
The coalgebra structure of $H^{*\rm cop}_{\mathcal D}$ is determined by
$\t(\xi)=\xi\ot\varepsilon+ \chi\ot \xi$,
$\t(\g)=\g\ot \g$ if $\mathcal{D}$ is of nilpotent type, and
$$\t(\g)=\g\ot \g +(\g^n(a)-1)\sum_{l+r=n}\frac{1}{(l)_{\rho}!(r)_{\rho}!}\g\chi^l\xi^r\ot\g\xi^l$$
if $\mathcal{D}$ is of non-nilpotent type, where $\g\in\G$.

\begin{proposition}\cite[Proposition 5]{KropRad}\label{2.1}
The Drinfeld double $D(H_{\mathcal D})$ is generated as an algebra by $G$, $x$, $\Gamma$ and $\xi$ subject to the relations defining $H_{\mathcal D}$ and $H^{*\rm cop}_{\mathcal D}$ and the following relations:
\begin{enumerate}
\item[(a)] $g\g=\g g$ for all $g\in G$ and $\g\in \Gamma$;
\item[(b)] $\xi g=\chi^{-1}(g)g\xi$ for all $g\in G$;
\item[(c)] $[x,\xi]=a-\chi$;
\item[(d)] $\g(a)x\g=\g x$ if $\mathcal D$ is {\bf nilpotent};
\item[(e)] $\g(a)x\g=\g x+\frac{\g^n(a)-1}{(n-1)_\rho!}\g(\rho a-\chi)\xi^{n-1}$  if $\mathcal D$ is {\bf non-nilpotent}.
\end{enumerate}
\end{proposition}

\section{Simple modules and projective modules}\label{s3}

Throughout this and the following sections, let $l(V)$ and ${\rm rl}(V)$ denote the length and the radical length (Loewy length) of a module $V$, respectively. For any integer $s\>0$, denote by $sV$ the direct sum of $s$ copies of $V$ (with the convention that 0$V$=0).

 For a module $V$, set $$V^x:=\{v\in V|xv=0\} \ {\rm and} \ V^{\xi}:=\{v\in V|\xi v=0\}.$$ Given elements $v_1, \cdots , v_s\in V$, we write $\langle v_1, \cdots, v_s\rangle$ for the submodule of $V$ generated by $\{v_1, \cdots, v_s\}$. Finally, for any set $X$, denote by $\sharp X$ the cardinality of $X$.

\subsection{Simple modules $V(l,\l)$}\label{s3.1}

~

Krop and Radford \cite[Subsection 2.2]{KropRad}, and independently Sun, Chen and Zhang \cite{SUN C Z}, classified the simple modules over $D(H_{\mathcal D})$ for any group datum $\mathcal{D}=(G,\chi,a,\a)$ with $G$ being abelian. In this subsection, we recall construction of the simple $D(H_{\mathcal D})$-modules.

Let $\Lambda:=\widehat{G\times \Gamma}$ be the set of characters of $G\times \Gamma$. Define a map ${\rm ev}_{a\chi^{-1}}: \L\rightarrow\Bbbk^{\times}$ by $${\rm ev}_{a\chi^{-1}}(\l)=\l(a\chi^{-1}), \ \l\in\L,$$ and let $K:={\rm Ker}({\rm ev}_{a\chi^{-1}})$.

Denote by $\widetilde{\L}$ the subset of characters $\l\in \L$ satisfying  $\l(a\chi^{-1})=\rho^s$ for some  $0\<s\<n-2$. For each $\l\in \L$, define $d(\l)=s$ if $\l\in \widetilde{\L}$, and   $d(\l)=-1$ otherwise. Let $\phi\in \L$ be defined by $\phi(g\g)=\chi^{-1}(g)\g(a), g\in G,\g\in \Gamma$, and define a mapping $\sigma: \Lambda \rightarrow \Lambda$ by $$\s(\l)=\l \phi^{d(\l)+1},$$ and set $\tau=\sigma^2$.

For any $1\<l\<n-1$, let $I_l=\{\l\in\widetilde{\L}|\l(a\chi^{-1})=\rho^{l-1}\}$. Define $$I_n:=\L\backslash\widetilde{\L}=\{\l\in\L|\l\notin\widetilde{\L}\},$$ and further decompose it as  $$I'_{n}=\{\l\in I_n|\l(a\chi^{-1})\neq \rho^k \ {\rm for \ any} \ 0\<k\<n-1\}$$ and $$I''_{n}=\{\l\in I_n|\l(a\chi^{-1})=\rho^{n-1}\}.$$ Then $I_n=I'_n\cup I''_n$, $\cup_{l=1}^{n-1} I_l=\widetilde{\L}$ and $\cup_{l=1}^{n} I_l={\L}$. By \cite[Theorem 2]{KropRad}, it follows that $\sharp I_l=\sharp I''_n=\sharp K$, where  $1\<l<n$.

For any $1\<l\<n$ and $\l\in I_l$, define $$\a_i(l,\l):=(i)_{\rho}(\l(\chi)-\l(a)\rho^{1-i})\in\Bbbk , i\>1.$$ Since $\s^2(\l)=\l\phi^n$, $\a_i(l,\tau(\l))=\a_i(l,\l)$. Now define $\b(l, \l)\in\Bbbk$ by
$$\b(1, \l):=1, \ \b(l, \l):=\a_1(l,\l)\a_2(l,\l)\cdots\a_{l-1}(l,\l) \ {\rm for} \  1<l\<n.$$
Note that $\a_i(l,\l)\neq 0$ for all $1\<i\<l-1$, and hence $\b(l,\l)\neq 0$.

Let $1\<l\<n$ and $\l\in I_l$. Let $\{v_i|0\<i\<l-1\}$ and $\{m_i|0\<i\<l-1\}$ be the natural and standard $\Bbbk$-bases of $V(l,\l)$ given in \cite{SUN C Z}, respectively.

 Let $1\<l\<n-1$. The $D(H_{\mathcal D})$-module action on $V(l,\l)$ is determined by
\begin{equation*}
\begin{array}{ll}
\vspace{0.2cm}
(g\g)v_i=(\phi^i\l)(g\g)v_i, & g\in G, \ \g\in\G, \ 0\<i\<l-1,\\
xv_i=\left\{
\begin{array}{ll}
v_{i+1},& 0\<i\<l-2,\\
0,& i=l-1,\\
\end{array}\right.&
\xi v_i=\left\{
\begin{array}{ll}
0, & i=0,\\
\a_i(l,\l)v_{i-1},& 1\<i\<l-1,\\
\end{array}\right.\\
\end{array}
\end{equation*}
or determined by
\begin{equation*}
\begin{array}{ll}
\vspace{0.2cm}
(g\g)m_i=(\phi^i\l)(g\g)m_i, & g\in G, \ \g\in\G, \ 0\<i\<l-1,\\
xm_i=\left\{
\begin{array}{ll}
\a_{i+1}(l,\l)m_{i+1},& 0\<i\<l-2,\\
0,& i=l-1,\\
\end{array}\right.&
\xi m_i=\left\{
\begin{array}{ll}
0, & i=0,\\
m_{i-1},& 1\<i\<l-1.\\
\end{array}\right.\\
\end{array}
\end{equation*}
Clearly,  $V(l,\l)^x=\Bbbk v_{l-1}=\Bbbk m_{l-1}$ and $V(l,\l)^{\xi}=\Bbbk v_0=\Bbbk m_0$.

Now we consider the case that $l=n$ and $\l\in I_n$. In this case, $V(n,\l)$ is projective. To describe the $D(H_{\mathcal D})$-module structure on $V(n,\l)$, it is necessary to distinguish between two cases for $\mathcal {D}$: nilpotent case and non-nilpotent case.

If $\mathcal{D}$ is of {\bf nilpotent type},  then the module action on $V(n,\l)$ is given by
\begin{equation*}
	\begin{array}{ll}
		\vspace{0.2cm}
		(g\g)v_i=(\phi^i\l)(g\g)v_i, & g\in G, \ \g\in\G, \ 0\<i\<n-1,\\
		xv_i=\left\{
		\begin{array}{ll}
			v_{i+1},& 0\<i\<n-2,\\
			0,& i=n-1,\\
		\end{array}\right.&
		\xi v_i=\left\{
		\begin{array}{ll}
			0, & i=0,\\
			\a_i(l,\l)v_{i-1},& 1\<i\<n-1,\\
		\end{array}\right.\\
	\end{array}
\end{equation*}
or given by
\begin{equation*}
\begin{array}{ll}
\vspace{0.2cm}
(g\g)m_i=(\phi^i\l)(g\g)m_i, &g\in G, \ \g\in\G, \ 0\<i\<n-1,\\
xm_i=\left\{
\begin{array}{ll}
\a_{i+1}(n,\l)m_{i+1},& 0\<i\<n-2,\\
0,& i=n-1,\\
\end{array}\right.&
\xi m_i=\left\{
\begin{array}{ll}
0, & i=0,\\
m_{i-1},& 1\<i\<n-1.\\
\end{array}\right.\\
\end{array}
\end{equation*}
Clearly,  $V(n,\l)^x=\Bbbk v_{n-1}=\Bbbk m_{n-1}$ and $V(n,\l)^{\xi}=\Bbbk v_0=\Bbbk m_0$.

If $\mathcal{D}$ is of {\bf non-nilpotent type}, then the module action on $V(n,\l)$ is given by
\begin{equation*}
	\begin{array}{ll}
		\vspace{0.2cm}
		(g\g)v_i=(\phi^i\l)(g\g)v_i, & g\in G, \ \g\in\G, \ 0\<i\<n-1,\\
		xv_i=\left\{
		\begin{array}{ll}
			v_{i+1},& 0\<i\<n-2,\\
			(\l^n(a)-1)v_0,& i=n-1,\\
		\end{array}\right.&
		\xi v_i=\left\{
		\begin{array}{ll}
			0, & i=0,\\
			\a_i(n,\l)v_{i-1},& 1\<i\<n-1,\\
		\end{array}\right.\\
	\end{array}
\end{equation*}
or given by
\begin{equation*}
\begin{array}{ll}
\vspace{0.2cm}
(g\g)m_i=(\phi^i\l)(g\g)m_i, &g\in G, \ \g\in\G, \ 0\<i\<n-1,\\
xm_i=\left\{
\begin{array}{ll}
\a_{i+1}(n,\l)m_{i+1},& 0\<k\<n-2,\\
\frac{\l^n(a)-1}{\b(n, \l)}m_0,& i=n-1,\\
\end{array}\right.&
\xi m_i=\left\{
\begin{array}{ll}
0, & i=0,\\
m_{i-1},& 1\<i\<n-1.\\
\end{array}\right.\\
\end{array}
\end{equation*}
Clearly, $V(n,\l)^{\xi}=\Bbbk v_0=\Bbbk m_0$.

\begin{proposition}\cite[Theorem 3.6]{SUN C Z}\label{3.6}
The following set $$\{V(l,\l)|1\<l\<n,\l\in I_l\}$$ is  a representative set of isomorphism classes of simple $D(H_{\mathcal D})$-modules. Moreover, $V(n,\l)$ is a projective $D(H_{\mathcal D})$-module for any $\l \in  I_n$.
\end{proposition}

\subsection{Indecomposable projective modules $P(l,\l)$}\label{s3.2}
~

Krop and Radford described all projective indecomposable $D(H_{\mathcal D})$-modules for any group datum $\mathcal{D}=(G,\chi,a,\mu)$ with $G$ being abelian in \cite[Subsection 2.3]{KropRad}. Sun, Chen and Zhang \cite{SUN C Z} reconstructed all indecomposable projective modules with matrices, as follows.

{\bf Case 1}: $\mathcal{D}$ is of {\bf nilpotent type}.
Let $1\<l\<n-1$, $\l \in I_l$, and let $\{v_0, v_1, \cdots, v_{n-1},\\ u_0, u_1, \cdots, u_{n-1}\}$ be the basis of $P(l,\l)$. The  $D(H_{\mathcal D})$-module action on $P(l,\l)$ is given by
\begin{equation*}
\begin{array}{ll}
\vspace{0.2 cm}
(g\g)v_i=(\phi^i\l)(g\g)v_i, & g\in G, \ \g\in\G, \ 0\<i\<n-1,\\
\vspace{0.2 cm}
(g\g)u_i=(\phi^{i-n+l}\l)(g\g)u_i, & g\in G, \ \g\in\G, \ 0\<i\<n-1,\\
\vspace{0.2 cm}
xv_i=\left\{
\begin{array}{ll}
v_{i+1}, & 0\<i\<n-2,\\
0,& i=n-1,\\
\end{array}\right.&
x u_i=\left\{
\begin{array}{ll}
u_{i+1}, & 0\<i\<n-2,\\
0,& i=n-1,\\
\end{array}\right.\\
\end{array}
\end{equation*}
\begin{equation*}
\begin{array}{ll}
\vspace{0.2 cm}
\xi v_i=\left\{
\begin{array}{ll}
u_{n-l-1},& i=0,\\
\a_i(l,\l)v_{i-1}+u_{n-l+i-1},& 1\<i\<l-1,\\
u_{n-1},& i=l,\\
\a_{i-l}(n-l,\s(\l))v_{i-1},& l+1\<i\<n-1,\\
\end{array}\right.&\\
\vspace{0.2 cm}
\xi u_i=\left\{
\begin{array}{ll}
0,& i=0,\\
\a_i(n-l,\s^{-1}(\l))u_{i-1}, & 1\<i\<n-l-1,\\
0,& i=n-l,\\
\a_{i-n+l}(l,\l)u_{i-1}, & n-l+1\<i\<n-1.\\
\end{array}\right.\\
\end{array}
\end{equation*}
Clearly, $P(l,\l)^{x}=\Bbbk v_{n-1}+\Bbbk u_{n-1}$ and $P(l,\l)^{\xi}=\Bbbk u_0+\Bbbk u_{n-l}$.

{\bf Case 2}: $\mathcal{D}$ is of {\bf non-nilpotent type}.
Let $1\<l\<n-1$, $\l \in I_l$, and let $\{v_0, v_1, \cdots, v_{n-1}, u_0, u_1, \cdots, u_{n-1}\}$ be the basis of  $P(l,\l)$. The $D(H_{\mathcal D})$-module action on $P(l, \l)$ is given by
\begin{equation*}
\begin{array}{ll}
\vspace{0.2 cm}
(g\g)v_i=(\phi^{i-n+l}\l)(g\g)v_i, & g\in G, \ \g\in\G, \ 0\<i\<n-1,\\
(g\g)u_i=(\phi^i\l)(g\g)u_i, & g\in G, \ \g\in\G, \ 0\<i\<n-1,\\
\end{array}
\end{equation*}
\begin{equation*}
\begin{array}{ll}
\vspace{0.2 cm}
x v_i=\left\{
\begin{array}{ll}
\a_{i+1}(n-l,\s^{-1}(\l))v_{i+1},& 0\<i\<n-l-2,\\
u_0,& i=n-l-1,\\
\a_{i+1-n+l}(l,\l)v_{i+1}+u_{i+1-n+l},& n-l\<i\<n-2,\\
y_{l,\l}v_0+u_l,& i=n-1,\\
\end{array}\right.&\\
x u_i=\left\{
\begin{array}{ll}
\a_{i+1}(l,\l)u_{i+1},& 0\<i\<l-2,\\
0,& i=l-1,\\
\a_{i+1-l}(n-l,\s(\l))u_{i+1}, & l\<i\<n-2,\\
z_{l,\l}u_0,& i=n-1,\\
\end{array}\right.\\
\end{array}
\end{equation*}
\begin{equation*}
\begin{array}{ll}
\xi v_i=\left\{
\begin{array}{ll}
0, & i=0,\\
v_{i-1},& 1\<i\<n-1,\\
\end{array}\right.&
\xi u_i=\left\{
\begin{array}{ll}
0, & i=0,\\
u_{i-1},& 1\<i\<n-1,\\
\end{array}\right.\\
\end{array}
\end{equation*}
where $y_{l,\l}=\frac{\rho^{1-l}\l(a)-\rho^l\l(\chi)}{(n-1)!_{\rho}}$ and $z_{l,\l}=\frac{\rho\l(a)-\l(\chi)}{(n-1)!_{\rho}}$. Moreover, $y_{l,\l}+z_{l,\l}=0$ and $z_{n-l,\s^{-1}(\l)}=y_{l,\l}$ by $\l(\chi)=\l(a)\rho^{1-l}$. Clearly, $P(l,\l)^x=\Bbbk u_{l-1}+\Bbbk(u_{n-1}-z_{l,\l}v_{n-l-1})$ and $P(l,\l)^{\xi}=\Bbbk v_0+\Bbbk u_0$.

\begin{proposition}\cite[Corollary 3.12]{SUN C Z}\label{3.12}
A representative set of isomorphism classes of indecomposable projective $D(H_{\mathcal D})$-modules is given by
 $$\{P(l,\l),V(n,\mu)|1\<l\<n-1, \l\in I_l, \mu \in I_n \}.$$
\end{proposition}

\section{The socle of the tensor product of two simple modules}\label{s400}

In this section, we investigate the socle of the tensor product of two simple non-projective modules.

For any integer $t$, let $c(t):=\left[\frac{t+1}{2}\right]$ denote the integer part of $\frac{t+1}{2}$. That is, $c(t)$ is the largest integer satisfying $c(t)\<\frac{t+1}{2}$.

\begin{lemma}\label{5.1}
Let $1\<t\<n$. Then  $\xi x^t=x^t\xi+(t)_\rho x^{t-1}(\chi-\rho^{1-t}a).$
\end{lemma}
\begin{proof}
We proceed by induction on $t$. If $t=1$, the  statement is obvious.
Assume $t>1$ and  that the lemma holds for $t-1$. Then
$$\begin{array}{rl}
\xi x^t=&x\xi x^{t-1}+\chi x^{t-1}-ax^{t-1}\\
=&x^t\xi+(t-1)_\rho x^{t-1}\chi-(t-1)_\rho \rho^{2-t}x^{t-1}a+\chi x^{t-1}-ax^{t-1}\\
=&x^t\xi+(t-1)_\rho x^{t-1}\chi+\rho^{t-1}x^{t-1}\chi-((t-1)_\rho\rho^{2-t}x^{t-1}a+\rho^{1-t}x^{t-1}a)\\
=&x^t\xi +(t)_\rho x^{t-1}\chi -\rho^{1-t}(t)_\rho x^{t-1}a\\
=&x^t\xi +(t)_\rho x^{t-1}(\chi-\rho^{1-t}a).
\end{array}$$
\end{proof}
For any $D(H_{\mathcal{D}})$-module $M$ and $\l\in \L$, set
$$M_{\l}:=\{v\in M|(g\g)v=\l(g\g)v, g\g\in G\times \G\}.$$
Then $M_{\l}$ is a subspace of $M$.
\begin{lemma}\label{5.2}
Let $M$ be a $D(H_{\mathcal{D}})$-module.
\begin{enumerate}
 \item[(1)]  Suppose there exists a $0\ne v\in M^{\xi}\cap M_{\l}$ for some $\l \in I_n$ such that $\xi(x^{k}v)\neq0$ for all $1\<k\<n-1$. Then $\langle v\rangle={\rm span}\{v, xv, \cdots, x^{n-1}v\}\cong V(n,\l)$.
 \item[(2)] Suppose there exists a $0\ne v\in M^{\xi}\cap M_{\l}$ for some $\l\in I_l$ with $1\<l\<n-1$. Then $\langle v\rangle$ is simple if and only if, for any $1\<k<n$, $x^kv\neq0$ implies $\xi(x^kv)\neq0$. In this case, $\langle v\rangle={\rm span}\{v, xv, \cdots, x^{l-1}v\}\cong V(l,\l)$.
\end{enumerate}
\end{lemma}
\begin{proof}
(1) By Lemma \ref{5.1},
$$\xi(x^{k}v) =(k)_{\rho}(\l(\chi)-\rho^{1-k}\l(a))x^{k-1}v=\a_k(n,\l)x^{k-1}v$$
for all $1\< k\< n-1$. We claim that $\g(x^kv)=\g^k(a)\l(\g)x^kv$ for all $\g\in\G$ and $0\<k\<n-1$.

If $\mathcal{D}$ is of nilpotent type, the claim follows from Proposition \ref{2.1}(d).  Assume now that  $\mathcal{D}$ is of non-nilpotent type.  Since $v\in M_{\l}$, we have $\g v=\l(\g)v$. Let $1\<k\<n-1$ and assume $\g(x^{k-1}v)=\g^{k-1}(a)\l(\g)x^{k-1}v$. Then
$$\g(x^{k}v) = \g(a)x\g(x^{k-1}v)-\frac{\g^n(a)-1}{(n-1)!_\rho}\g (\rho a-\chi)\xi^{n-1}(x^{k-1}v)=\g^k(a)\l(\g)(x^kv).$$
Thus,  the claim holds.

 On the other hand,  since $v\in M_{\l}$ and $gx=\chi^{-1}(g)xg$ for $g\in G$, we obtain
 $$g(x^kv)=\chi^{-k}(g)\l(g)x^kv\ \text{for all} \ 0\<k\<n-1.$$
 Therefore,
 $$(g\g)(x^kv)=(\phi^k\l)(g\g)x^kv\ \text{for all}\  0\<k\<n-1, \ g\g\in G\times\G.$$
 If $\mathcal D$ is of nilpotent type, then $x^n=0$, and so $x^nv=0$. \\
 If $\mathcal D$ is of non-nilpotent type, then $x^n=a^n-1$, hence,  $x^nv=(\l^n(a)-1)v$. It follows that
$$\langle v\rangle={\rm span}\{ v, xv,\cdots ,x^{n-1}v \}\cong V(n,\l).$$

(2) Let $0\ne v\in M^{\xi}\cap M_{\l}$ for some $\l\in I_l$ with $1\<l\<n-1$.
Assume first that $\langle v\rangle$ is simple. Then $\langle v\rangle^{\xi}=\Bbbk v$ by \cite[Corollary 3.7]{SUN C Z}. If $x^{k}v\neq 0$ for some $1\<k \< n-1$, then $x^{k}v\notin\langle v\rangle^{\xi}$
since $av=\l(a)v$ and $a(x^kv)=\rho^{-k}\l(a)x^kv$.
Hence,  $\xi(x^kv)\neq 0$.

Conversely, assume that  whenever $x^{k}v\neq 0$ for any $1\<k<n$, we have $\xi (x^{k}v)\ne 0$. Note that $x^nv=0$.  Let $t$
be the minimal positive integer such that $ x^tv= 0$. Then $1\<t\<n$.\\
If $t=1$, then $\langle v\rangle=\Bbbk v\cong V(1,\l)$. Now assume
$1<t\< n$. Then for any $1\< k\< t-1$, $x^{k}v\neq 0$ and hence $\xi(x^{k}v)\neq 0$. By  Lemma \ref{5.1}, we have
$$\xi (x^{k}v)=(k)_{\rho}(\l(\chi)-\rho^{1-k}\l(a))x^{k-1}v \ \text{for all}\ k\>1.$$
 In particular,
$$0 =\xi (x^{t}v)=(t)_{\rho}(\l(\chi)-\rho^{1-t}\l(a))x^{t-1}v.$$
Hence,  $(t)_{\rho}(\l(\chi)-\rho^{1-t}\l(a))=0$.\\
 If $t\neq l$ , then $\l(\chi)-\rho^{1-t}\l(a)\neq 0$  since $1\<l,t\< n$ and $\l(\chi)=\rho^{1-l}\l(a)$. Thus,  $(t)_{\rho}=0$, and so $t=n$, leading $1\<l<t$.  This implies $a^lv\neq 0$. But $\xi(x^{l}v)=0$, which contradicts the assumption. Therefore,  $t=l$ and
$$\xi (x^{k}v) = \alpha _{k}( l,\l )x^{k-1}v\ \text{for}\ 1\<k\<l-1.$$
As in (1),  we can show that $(g\g)(x^kv)=(\phi^k\l)(g\g)x^kv$ for all $0\<k\<l-1$ and $g\g\in G\times\G$. Hence,
$$\langle v\rangle={\rm span}\{ v, xv,\cdots ,x^{l-1}v \}\cong V(l,\l).$$
\end{proof}

\begin{lemma}{\label{5.3}}
Let $1\<l\<n$, $1\<l'\<n-1$, $\l\in I_1$, $\mu \in I_l$ and $\l'\in I_{l'}$. Then
$$V(1,\l)\ot V(l,\mu)\cong V(l,\l\mu) \ \text{ and }\ V(1,\l)\ot P(l',\l')\cong P(l,\l\l').$$
\end{lemma}
\begin{proof}
The assertion follows by a straightforward verification.
\end{proof}

In the remainder of this section, unless otherwise stated, assume that $2\< l\< l'\<n-1$, $\l\in I_l$ and $\l'\in I_{l'}$. Let
  $$M=V(l,\l)\otimes V(l',\l').$$
Let $\{m_0, \cdots, m_{l-1}\}$ and $\{m'_0, \cdots, m'_{l'-1}\}$ be the standard bases of $V(l, \l)$ and $V(l', \l')$, respectively.
Then $$\{m_i\ot m'_j\mid 0\<i\<l-1, 0\<j\<l'-1\}$$
 is a basis of $M$.  For each integer $s$ with $0\<s\<l+l'-2$, set
 $$M_{[s]}={\rm span}\{m_i\ot m'_j\mid i+j=s\}$$
  and define $M_{[-1]}=M_{[l+l'-1]}=0$.

\begin{lemma}\label{5.4}
Retain the above notations.
\begin{enumerate}
 \item[(1)] For any $0\<s\<l+l'-2$, we have $xM_{[s]}\subseteq M_{[s+1]}$ and $\xi M_{[s]}\subseteq M_{[s-1]}$.
 \item[(2)] The subspace $M^{x}$ decomposes as
 $$M^{x}=\oplus_{s=l'-1}^{l+l'-2}\left(M^{x}\cap M_{[s]}\right)$$ and ${\rm dim}(M^{x}\cap M_{[s]})=1$ for any $l'-1\<s\<l+l'-2$.
 \item[(3)] The subspace $M^{\xi}$ decomposes as
 $$M^{\xi}=\oplus _{s=0}^{l-1}\left(M^{\xi}\cap M_{[s]}\right)$$ and ${\rm dim}(M^{\xi}\cap M_{[s]})=1$ for any $0\<s\<l-1$.
 \item[(4)]  If $l<l'$, then $M^{x}\cap M^{\xi}=0$;  if $l= l'$,  then ${\rm dim}(M^{x}\cap M^{\xi})=1$.
\end{enumerate}
\end{lemma}

\begin{proof}
The statements follow by direct computation using the action of $x$ and $\xi$ on the basis elements of $M$.
\end{proof}

By Lemma \ref{5.4}, we have ${\rm dim}(M^{\xi}\cap M_{[s]})=1$ for all $0\<s\<l-1$. In the rest of this section,  let  $0 \neq x_s\in M^{\xi}\cap M_{[s]}$ and let $U_s=\langle x_s\rangle$  denote the submodule of $M$ generated by $x_s$, $0\<s\<l-1$.

\begin{lemma}\label{5.5}
Let $U$ be a simple submodule of $M$. Then $U=U_s$ for some $0\<s\<l-1$.
\end{lemma}

\begin{proof} The assertion follows directly from Lemma \ref{5.4} together with the known structure of simple $D(H_{\mathcal{D}})$-modules.
\end{proof}

\begin{lemma}\label{5.6}
Let  $0\< s\< l-1$. Then $x_s\in M_{\phi^s\l\l'}$ and  $(\phi^s\l\l')(a\chi^{-1})=\rho^{l+l'-2-2s}$.
\end{lemma}

\begin{proof}
Let $0\<i\<s$ and $g\g\in G\times\G$. Then
 \begin{eqnarray*}
 g (m_i\ot m'_{s-i}) &=& g m_i\ot gm'_{s-i}\\
 &=& (\phi^i\l)(g)m_i\ot(\phi^{s-i}\l')(g)m'_{s-i}\\
 &=& (\phi^s\l\l')(g)m_i\ot m'_{s-i}.
 \end{eqnarray*}
 If $\mathcal{D}$ is of nilpotent type, then
 $$\g (m_i\ot m_{s-i})=\g m_i\ot \g m_{s-i}=(\phi^s\l\l')(\g)m_i\ot m'_{s-i}.$$
  Now assume that $\mathcal{D}$ is of non-nilpotent type. Then,  we have
$$\begin{array}{rl}
&\g (m_{i}\otimes m'_{s-i}) \\
=&\g m_{i}\otimes \g m'_{s-i}+(\g^{n}(a)-1)\sum_{r+l=n}\frac{1}{(l)!_\rho (r)!_\rho}(\g\chi^l\xi^rm_i\ot \g\xi^lm'_{s-i})\\
=&(\phi^s\l\l')(\g)m_{i}\otimes  m'_{s-i}+(\g^{n}(a)-1)\sum_{r+l=n}\frac{1}{(l)!_\rho (r)!_\rho}(\g\chi^l\xi^rm_i\ot \g\xi^lm'_{s-i}).\\
\end{array}$$
If $\g\chi^l\xi^rm_i\ot \g\xi^lm'_{s-i}\neq 0$, then  $i\>r$ and $s-i\>l$, hence $s\>n$ since $r+l=n$, a contradiction.  Therefore,  all such terms vanish, and thus
$$\g (m_i\ot m'_{s-i})=(\phi^s\l\l')(\g) m_i\ot m'_{s-i}.$$
It follows that $(g\g)x_s=(\phi^s\l\l')(g\g)x_s$;  hence $x_s\in M_{\phi^s\l\l'}$.
Finally,  by direct computation,  $(\phi^s\l\l')(a\chi^{-1})=\rho^{l+l'-2-2s}$. This completes the proof.
\end{proof}

\begin{lemma}\label{5.7}
Let $t=l+l'-n-1$.
\begin{enumerate}
\item[(1)] Suppose $t< 0$ and let $0\< s\< l-1$. Then $1\<l+l'-1-2s\<n-1$, and \\
$U_s\cong V(l+l'-1-2s, \phi^s\l\l')$.
\item[(2)] Suppose $t\>0$ is even.
\begin{enumerate}
\item[(a)] If $c(t)+1\< s \< l-1$, then $1\<l+l'-1-2s\<n-2$,  and \\
$U_s\cong V(l+l'-1-2s, \phi^s\l\l')$.
\item[(b)] If $s=c(t)$, then $l+l'-1-2s=n$, and
$U_s\cong V(n, \phi^s\l\l')$.
\item[(c)]  If $0\<s\<c(t)-1$, then $U_s$ is not simple.
\end{enumerate}
\item[(3)] Assume that $t\>0$ is odd.
\begin{enumerate}
\item[(a)]  If $c(t)\< s \<l-1$, then $1\<l+l'-1-2s\<n-1$, and\\
 $U_s\cong V(l+l'-1-2s, \phi^s\l\l')$.
\item[(b)] If $0\<s\<c(t)-1$, then $U_s$ is not simple.
\end{enumerate}
\end{enumerate}
\end{lemma}

\begin{proof}
(1) Since $t<0$ and $0\<s\<l-1$,  we have
$$1\<l'-l+1\<l+l'-1-2s\<l+l'-1\<n-1.$$
 Hence, by Lemma \ref{5.6}, we obtain  $\phi^s\l\l'\in I_{l+l'-1-2s}$.
 Then by Lemmas \ref{5.1} and \ref{5.6},  for all $k\>1$,
 $$\xi(x^kx_s)=\a_k(l+l'-1-2s, \phi^s\l\l')x^{k-1}x_s.$$
By induction on $k$, it follows that
$$\xi(x^kx_s)\neq 0 \ \text{and}\ x^kx_s\neq 0 \ \text{for all} \ 1\<k\<l+l'-2-2s.$$
However, $\xi(x^{l+l'-1-2s}x_s)=0$, and hence $x^{l+l'-1-2s}x_s\in M^{\xi}$.

By Lemma \ref{5.4}(1), this implies that $x^{l+l'-1-2s}x_s\in M^{\xi}\cap M_{[l+l'-1-s]}$.
Since $0\<s\<l-1$, we have $l-1<l+l'-1-s\<l+l'-1$.  Therefore,  by Lemma \ref{5.4}(3), $M^{\xi}\cap M_{[l+l'-1-s]}=0$.  It follows  that $x^{l+l'-1-2s}x_s=0$.

Applying Lemma \ref{5.2}(2),  we conclude that $U_s\cong V(l+l'-1-2s, \phi^s\l\l')$.

(2) Suppose first $c(t)+1\<s\<l-1$. Since $t\>0$ is even, we have  $1\<l+l'-1-2s\<l+l'-1-2(c(t)+1)=n-2$ and $l-1<l+l'-1-s\<l+l'-2-c(t)\<l+l'-2$. Then, by  an argument identical  to that of part (1), we obtain $U_s\cong V(l+l'-1-2s, \phi^s\l\l')$.

Next, suppose $s=c(t)$. Then $l+l'-1-2s=n$, and hence,  by Lemma \ref{5.6}, $\phi^s\l\l'\in I_n$. Repeating the reasoning of part (1), we find that
$$x^kx_s\neq 0\ \text{and}\ \xi(x^kx_s)\neq 0\  \text{for all}  1\<k<n. $$
Therefore,  by Lemma \ref{5.2}(1),  $U_s\cong V(n, \phi^s\l\l')$.

Finally, suppose $0\<s\<c(t)-1$. Then
$$n+2\<l+l'-1-2s\<l+l'-1\<2n-3,$$
so that $2\<l+l'-1-n-2s\<n-3$. By Lemma \ref{5.6}, it follows that $\phi^s\l\l'\in I_{l+l'-1-n-2s}$.

As in part (1), one checks that
$$x^kx_s\neq 0\ \text{and}\ \xi(x^kx_s)\neq 0\ \text{for all}\ 1\<k\<l+l'-2-n-2s,$$
while $\xi(x^{l+l'-1-n-2s}x_s)=0$. In particular, $x^{l+l'-2-n-2s}x_s\neq 0$. By Lemma \ref{5.4}(1), this element lies in $M_{[l+l'-2-n-s]}$.
Since $l+l'-2-n-s\<l+l'-2-n\<l'-3$, we have $x^{l+l'-2-n-2s}x_s\notin M^x$ by Lemma \ref{5.4}(2), and  hence $x^{l+l'-1-n-2s}x_s\neq 0$.  It follows from Lemma \ref{5.2}(2) that $U_s$ is not simple.

(3) When $t>0$ is odd, the proof is entirely analogous to that of part (2):  if
$c(t)\< s \<l-1$, then  $$U_s\cong V(l+l'-1-2s, \phi^s\l\l'),$$
while for  $0\<s\<c(t)-1$,  the submodule $U_s$ is not simple.
\end{proof}

%\begin{corollary}{\label{40.8}}
%Let  $t\> 0$ and $2\< s\< c(t)+1$. Then the submodule $\langle x_s\rangle$ of $M$ is not %simple.
%\end{corollary}
%\begin{proof}
%  It follows from Lemmas \ref{5.2}, \ref{5.5} and \ref{5.7}.
%\end{proof}
%\begin{lemma}{\label{40.9}} Let $2\<s\<l+1$ and $1\<k\<n-1$. Then
%\item[(1)]$g(x^{k}x_{s})=\chi^{2-s-k}(g)\l(g)\l'(g)x^{k}x_{s}$.
%\item[(2)]$\g(x^{k}x_{s})=\g^{s+k-2}(a)\l(\g)\l'(\g)x^{k}x_{s}$.
%\end{lemma}
%\begin{proof}
%(1) It follows from Lemma \ref{5.6} and a  straightforward computation.
%(2) If $\mathcal{D}$ is nilpotent, then (2) follows from Lemma \ref{5.6}. Now assume that $\mathcal{D}$ is non-nilpotent. Then by the algebra structure of $D(H_{\mathcal{D}})$, we have
%$$\g(x^{k}x_{s}) = \g(a) x\g x^{k-1}x_{s}-\frac{\g^n(a)-1}{(n-1)!_\rho}\g(\rho a-\chi)\xi^{n-1}x^{k-1}x_{s}.$$
%By Lemmas \ref{5.1} and \ref{5.6}, $\xi (x^{k}x_{s})=(k)_{\rho}(\chi^{2-s}\widehat{a}^{s-2}\l\l'(\chi)-\rho^{1-k}\chi^{2-s}\widehat{a}^{s-2}\l\l'(a))x^{k-1}x_{s}$. Thus, $\xi^{n-1}$ $x^{k-1}x_{s}=0$ by $1\<k\<n-1$. Consequently, it follows from Lemma \ref{5.6} that $\g(x^{k}x_{s}) = \g^{k}(a)x^{k}\g x_{s}=\g^{s+k-2}(a)\l(\g)\l'(\g)x^{k}x_{s}$.
%\end{proof}

\begin{corollary}\label{5.8}
Let $t=l+l'-n-1$.
\begin{enumerate}
\item[(1)]If $t<0$, then ${\rm soc}M \cong  \oplus _{s=0}^{l-1 } V(l+l^{\prime }-1-2s,\phi^{s}\l\l' ).$
\item[(2)]If $t\>0$, then
${\rm soc}M\cong\oplus_{s=c(t)}^{l-1 } V(l+l^{\prime }-1-2s,\phi^{s}\l\l' )$.
\end{enumerate}
\end{corollary}

\begin{proof}
It follows from Lemma \ref{5.6} and Lemma \ref{5.7}.
\end{proof}

\section{ The decomposition rules for tensor product modules}\label{s500}

In this section, we study the tensor products of simple  and indecomposable projective  $D(H_{\mathcal{D}})$-modules, and determine their decompositions into the direct sums of indecomposable modules.

\subsection{Tensor products of two simple modules}\label{s50.1}
~

In this subsection, we investigate the tensor products of two simple modules. We first consider the tensor products of two non-projective simple modules.  Note that  for any $D(H_{\mathcal D})$-modules, there is an isomorphism between $M\ot N$ and $N\ot M$.

Throughout the following, we write  $P(n,\l):=V(n,\l)$ for any $\l\in I_n$.

\begin{lemma}\label{6.1}
Let $1\<l, l'\<n-1$, $\l\in I_l$, $\l'\in I_{l'}$ and $s\>1$. Then
$$V(l,\l)\ot \Omega^{\pm s}V(l',\l')\cong\Omega^{\pm s}(V(l,\l)\ot V(l',\l'))\oplus P$$ for some projective module $P$.
\end{lemma}

\begin{proof}
It follows from the fact that the tensor product of a projective module with any module is projective.
\end{proof}

\begin{theorem}\label{6.2}
Suppose that $1\<l\<l'<n$, $\l\in I_{l}$ and $\l'\in I_{l'}$.
Let $t=l+l'-n-1$.
\begin{enumerate}
\item[(1)] If  $t<0$, then $V(l,\l)\otimes V(l',\l')\cong \oplus _{s=0}^{l-1}V( l+l'-1-2s,\phi^s\l\l')$.
\item[(2)] If $t\>0$, then $$\begin{array}{rcl}
V(l,\l)\otimes V(l',\l')
&\cong &(\oplus_{s=c(t)}^tP(l+l'-1-2s,\phi^s\l\l'))\\
&&\oplus(\oplus_{t+1\<s\<l-1}V(l+l'-1-2s,\phi^s\l\l')).\end{array}$$
\end{enumerate}
\end{theorem}

\begin{proof}
(1) Assume $t<0$.   By Corollary \ref{5.8}, we have
$${\rm soc}(V(l,\l)\otimes V(l',\l'))\cong\oplus_{s=0}^{l-1}V(l+l'-1-2s, \phi^s\l\l').$$
Since both sides have the same dimension, it follows that
$$V(l,\l)\otimes V(l',\l')\cong\oplus_{s=0}^{l-1}V(l+l'-1-2s, \phi^s\l\l').$$

(2) Assume $t\>0$.  By \cite[Lemma 3.11]{SUN C Z}, there exists an exact sequence
$$0\rightarrow V(n-l', \s(\tau^k(\l')))\rightarrow\Omega V(n-l',\s(\tau^k(\l')))\rightarrow V(l',\tau^k(\l'))\oplus V(l',\tau^{k+1}(\l') )\rightarrow 0.$$
Taking the direct sum of these sequences over all $k$,  and then tensoring on the left by $V(l,\l)$, we obtain another exact sequence: \\
$$\begin{array}{rl}
0\rightarrow &\oplus_{k=0}^{m-1}V(l,\l)\ot V(n-l',\s(\tau^k(\l')))\\
\rightarrow &\oplus_{k=0}^{m-1}V(l,\l)\ot\Omega V(n-l',\s(\tau^k(\l')))\\
\rightarrow &2(\oplus_{k=0}^{m-1}V(l,\l)\ot V(l',\tau^k(\l')))
\rightarrow 0.
\end{array}$$

Since $t\>0$,  we have $n-l'<l$. Then by (1), we have
$$\begin{array}{rl}
&V(l,\l)\ot V(n-l',\s(\tau^k(\l')))\\
\cong &\oplus_{s=0}^{n-l'-1}V(n-l'+l-1-2s,\phi^{s}\s(\tau^k(\l'))\l)\\
\cong &\oplus_{s=l'}^{n-1}V(n+l'+l-1-2s,\phi^{kn+s}\l'\l ).
\end{array}$$
Hence, by Lemma \ref{6.1}, we obtain
$$\begin{array}{rl}
&\oplus_{k=0}^{m-1}V(l,\l)\ot \Omega V(n-l',\s(\tau^k(\l')))\\
\cong &(\oplus_{k=0}^{m-1} \oplus_{s=l'}^{n-1}\Omega V(n+l'+l-1-2s,\phi^{kn+s}\l'\l))\oplus P
\end{array}$$
for some projective module $P$. Let
$$L=\oplus_{k=0}^{m-1} \oplus_{s=l'}^{n-1} V(n+l'+l-1-2s,\phi^{kn+s}\l'\l),  M=\oplus_{k=0}^{m-1}V(l,\l)\ot V(l',\tau^k(\l')).$$
 Then,  the  exact sequence obtained above takes the form:
\begin{eqnarray*}
0\rightarrow L\xrightarrow{f} \Omega L\oplus P \xrightarrow{g} 2M\rightarrow 0.
\end{eqnarray*}
For $l'\<s\<n-1$, we have $2\<n+l'+l-1-2s\<n-1$. Hence, $L$ is  semisimple, and any simple submodule of $L$ is not projective. It follows from \cite[Lemma 3.11]{SUN C Z} that $L\cong {\rm soc}(\Omega L)$. Moreover, by \cite[Proposition 3.2 (1)]{SUN C Z}, each simple submodule of $L$ appears with  multiplicity one as a composition factor of $L$.

By Corollary \ref{5.8}, we have
$$\text{soc} M\cong \oplus_{k=0}^{m-1}\oplus_{s=c(t)}^{l-1}V(l+l'-1-2s, \phi^{kn+s}\l\l').$$
According to \cite[Proposition 3.2 (1)]{SUN C Z}, $L$ and soc$M$ share no isomorphic simple submodule.

Let $\pi:$ $\Omega L\oplus P\rightarrow P$ be the canonical projection. If $\pi f\neq 0$, then there exists a simple submodule $S$ of $L$ such that $\pi (f(S))\neq 0$, implying $\pi (f(S))\cong S$. Consequently,  as a composition factor of soc$(\Omega L\oplus P)$, $S$ would appear with multiplicity at least $2$. This forces $S$ to be isomorphic to a simple submodule of $M$,  contradicting the above observation. Hence $\pi f=0$, and therefore ${\rm Im}(f)\subseteq \Omega L$. It follows that ${\rm Im}(f)={\rm soc}(\Omega L)$.

By \cite[Lemma 3.11]{SUN C Z}, we thus obtain
$$\begin{array}{rl}
2M &\cong \Omega L/{\rm soc}(\Omega L)\oplus P\\
&\cong P \oplus (\oplus_{k=0}^{m-1} \oplus_{s=l'}^{n-1}2 V(2s+1-l'-l,\phi^{kn+l+l'-1-s}\l\l'))\\
&\cong P \oplus (\oplus_{k=0}^{m-1} \oplus_{s=t+1}^{l-1}2V(l+l'-1-2s,\phi^{kn+s}\l\l')).
\end{array}$$
 Thus, there exists a projective module $Q$ such that
 $$M\cong Q\oplus (\oplus_{k=0}^{m-1} \oplus_{s=t+1}^{l-1}V(l+l'-1-2s,\phi^{kn+s}\l\l')).$$
 From the description of soc$M$ given above, we obtain
 $$\text{soc}Q\cong \oplus_{k=0}^{m-1} \oplus_{s=c(t)}^{t}V(l+l'-1-2s,\phi^{kn+s}\l\l').$$
 Since $Q$ is a projective cover (injective envelope ) of soc$Q$,  we have
  $$Q\cong \oplus_{k=0}^{m-1} \oplus_{s=c(t)}^{t}P(l+l'-1-2s,\phi^{kn+s}\l\l')$$
   Therefore
$$\begin{array}{rl}
  M\cong& (\oplus_{k=0}^{m-1}\oplus_{s=c(t)}^{t}P(l+l'-1-2s,\phi^{kn+s}\l\l'))\\
  &\oplus(\oplus_{k=0}^{m-1}
  \oplus_{s=t+1}^{l-1}V(l+l'-1-2s,\phi^{kn+s}\l\l')).\end{array}$$
  Since all the summands in $M$ are pair-wise non-isomorphic,  and
  $$\text{soc}(V(l,\l)\ot V(l',\l'))\cong \oplus_{s=c(t)}^{l-1}V(l+l'-1-2s,\phi^{s}\l\l'),$$
   we finally obtain
$$\begin{array}{rl}
V(l,r)\ot V(l',r')\cong &(\oplus_{s=c(t)}^t P(l+l'-1-2s,\phi^{s}\l\l'))\\
&\oplus(\oplus_{s=t+1}^{l-1}V(l+l'-1-2s,\phi^{s}\l\l')).
\end{array}$$
\end{proof}

 Next we consider the tensor products of non-projective simple modules with projective simple modules. Let $N=V(l,\l)\otimes V(n,\l')$, where $1\<l\<n-1$, $\l\in I_l$ and $\l'\in I_{n}$. \\
  Let $\{m_0, \cdots, m_{l-1}\}$ and $\{m'_0, \cdots, m'_{n-1}\} $ be the standard bases of $V(l,\l)$ and $V(n,\l')$, respectively. Then
  $$\{m_i\ot m'_j\mid 0\<i\<l-1, 0\<j\<n-1\}$$
  forms a basis of $N$.  For $0\<s\<l+n-2$, set $N_{[s]}={\rm span}\{m_i\ot m'_j\mid i+j=s\} $  and define $N_{[-1]}=N_{[l+n-1]}=0$.

If $\mathcal D$ is of non-nilpotent type and $\l'\in I''_n$, then $\l'(a)=\rho^{n-1}\l'(\chi)$, so that
$$\l'(a)^n=\rho^{(n-1)n}\l'(\chi)^n=\l'(\chi^n)=1.$$
 This implies that $xm'_{n-1}=0$. Therefore,  by an argument similar to Lemma \ref{5.4}, we obtain the following result.

\begin{lemma}\label{6.3}
Retain the above notations.
\begin{enumerate}
\item[(1)]$xN_{[s]}\subseteq N_{[s+1]}$ for $0\<s\<n-2$. Moreover, if $\l'\in I''_n$, then $xN_{[s]}\subseteq N_{[s+1]}$ for $0\<s\<n+l-2$
%\item[(2)] ${\rm dim}(N^x\cap N_{[s]})=0$ for all $0\<s\<n-2$.
\item[(2)] $\xi N_{[s]}\subseteq N_{[s-1]}$ for $0\<s\<l+n-2$.
\item[(3)]$N^{\xi}=\oplus _{k=0}^{l-1}\left( N^{\xi}\cap N_{[s]}\right)$, and ${\rm dim}(N^{\xi}\cap N_{[s]})=1$ for $0\<s\<l-1$.
\end{enumerate}
\end{lemma}

In what follows, choose nonzero element $x_s\in N^{\xi}\cap N_{[s]}$ for $0\<s\<l-1$, and let $U_s:=\langle x_s\rangle$  denote the submodule of $N$ generated by $x_s$.  Analogously to Lemma \ref{5.5}, we obtain the following lemma.

\begin{lemma}\label{6.4}
Let $U$ be a simple submodule of $N$. Then $U=U_s$ for some $0\<s\<l-1$.
\end{lemma}

\begin{lemma}\label{6.5}
Let $0\<s\<l-1$. Then $x_s\in N_{\phi^s\l\l'}$. Moreover,
\begin{enumerate}
\item[(i)] if $\l'\in I''_n$, then $(\phi^s\l\l')(a\chi^{-1})=\rho^{l-2-2s}$;
\item[(ii)]   if $\l'\in I'_n$, then  $\phi^s\l\l'\in I'_n$.
\end{enumerate}
\end{lemma}

\begin{proof}
Similar to the proof of Lemma \ref{5.6}.
\end{proof}

\begin{lemma}\label{6.6}
Let $0\<s\<l-1$.
\begin{enumerate}
\item[(1)] Suppose that $\l'\in I'_n$. Then $U_s\cong V(n, \phi^s\l\l')$.
\item[(2)] Suppose that $\l'\in I''_n$ and $l$ is even. If $\frac{l}{2} \<s\<l-1$, then $2\<n+l-1-2s\<n-1$ and $U_s\cong V(n+l-1-2s, \phi^s\l\l')$.
\item[(3)] Suppose that $\l'\in I''_n$ and $l$ is odd.  If $\frac{l+1}{2}\<s\<l-1$ , then $2\<n+l-1-2s\<n-2$ and $U_s\cong V(n+l-1-2s, \phi^s\l\l')$. If $s=\frac{l-1}{2}$, then  $U_s\cong V(n, \phi^s\l\l')$.
\end{enumerate}
\end{lemma}

\begin{proof}
By Lemmas \ref{5.1} and \ref{6.5}, we have
\begin{eqnarray}\eqlabel{eq6.1}
\xi(x^kx_s)=(k)_\rho((\phi^s\l\l')(\chi)-\rho^{1-k}(\phi^s\l\l')(a))x^{k-1}x_s,\ k\>1.
\end{eqnarray}

(1)  Since  $\l'\in I'_n$, it follows from  Lemma \ref{6.5} that $\phi^s\l\l'\in I'_n$. Hence
$$(\phi^s\l\l')(\chi)-\rho^{1-k}(\phi^s\l\l')(a)\neq 0\ \text{for all}\  k\in \mathbb{Z}.$$
By Equation \equref{eq6.1} and induction on $k$, we obtain that
$\xi(x^kx_s)\neq 0$ for all $1\<k\<n-1$.  Therefore,  $U_s\cong V(n, \phi^s\l\l')$ by Lemmas \ref{5.2}(1) and \ref{6.5}.

(2) Suppose $\frac{l}{2}\<s\<l-1$. Since $\l'\in I''_n$ and $l$ is even, we have
$$2\<n+l-1-2s\<n-1\ \text{and}\ \phi^s\l\l'\in I_{n+l-1-2s}$$
 by Lemma \ref{6.5}. Then by Eq.\equref{eq6.1}, we have
$$\xi(x^kx_s)=\a_k(n+l-1-2s, \phi^s\l\l')x^{k-1}x_s, k\>1.$$
Thus, by induction on $k$, we have $\xi(x^kx_s)\neq 0$ for all $1\<k\<n+l-2-2s$, while  $\xi(x^{n+l-1-2s}x_s)=0$.  Since $n\<n+l-1-s\<n+l-2$, it follows from Lemma \ref{6.3}(1) and (3) that  $x^{n+l-1-2s}x_s\in N^{\xi}\cap N_{[n+l-1-s]}=0$, and hence $x^{n+l-1-2s}x_s=0$.

 Consequently,  $U_s\cong V(n+l-1-2s, \phi^s\l\l')$ by Lemmas \ref{5.2}(2) and \ref{6.5}.

(3) The proof  is similar to that of $(2)$.
This  complete the proof.
\end{proof}

\begin{theorem}\label{6.7}
Retain the above notations.
\begin{enumerate}
\item[(1)] If $\l'\in I'_n$, then $N\cong\oplus_{s=0}^{l-1} V(n,\phi^s\l\l')$.
\item[(2)] If $\l'\in I''_n$, then $N\cong\oplus _{s=c(l-1) }^{l-1} P(n+l-1-2s,\phi^{s}\l\l')$.
    \end{enumerate}
\end{theorem}

\begin{proof}
(1) If $\l'\in I'_n$, then by Lemmas \ref{6.4} and \ref{6.6}(1), we have
$${\rm soc}N=\oplus_{s=0}^{l-1}U_s\cong\oplus_{s=0}^{l-1}V(n,\phi^s\l\l').$$
Hence, ${\rm dim}({\rm soc}N)=ln={\rm dim}N$, which implies $N={\rm soc}N\cong\oplus_{s=0}^{l-1} V(n,\phi^s\l\l')$.

(2) Assume $\l'\in I''_n$. We only consider the case when $l$ is odd,  since the even case is similar. By Lemma \ref{6.6}(3),
$$\oplus_{s=\frac{l-1}{2}}^{l-1} U_s\cong\oplus_{s=\frac{l-1}{2}}^{l-1}V(n+l-1-2s,\phi^s\l\l').$$
Since $\oplus_{s=\frac{l-1}{2}}^{l-1}U_s\subseteq{\rm soc}N$ and $N$ is injective, it follows that $N$ contains a submodule $L\cong\oplus_{s=\frac{l-1}{2}}^{l-1}P(n+l-1-2s,\phi^s\l\l')$.  Comparing dimensions,  ${\rm dim}L=ln={\rm dim}N$, we conclude
$$N=L\cong\oplus_{s=\frac{l-1}{2}}^{l-1}P(n+l-1-2s,\phi^s\l\l').$$
\end{proof}

Next, we consider the tensor products of two simple projective modules.
Let $T=V(n,\l)\otimes V(n,\l')$, where $\l, \l'\in I_n$. Let $\{m_0, \cdots, m_{n-1}\}$ and $\{m'_0, \cdots, m'_{n-1}\}$ denote the standard bases of $V(n,\l)$ and $V(n,\l')$, respectively. Then $\{m_i\ot m'_j\mid 0\<i, j\<n-1\}$ is a basis of $T$.

 For $0\<s\<2n-2$, set
 $$T_{[s]}={\rm span}\{m_i\ot m'_j\mid i+j=s\}\ \text{and}\ T_{[-1]}=T_{[2n-1]}=0.$$

Similarly to Lemma \ref{5.4}, we have the following result.

\begin{lemma}\label{6.8}
Retain the above notations.
\begin{enumerate}
\item[(1)] $xT_{[s]}\subseteq T_{[s+1]}$ for $0\<s\<n-2$, and $xT_{[s]}\subseteq T_{[s+1]}\oplus T_{[s+1-n]}$ for $n-1\<s\<2n-2$.
\item[(2)] If $\mathcal D$ is of nilpotent type, then $xT_{[s]}\subseteq T_{[s+1]}$ for all $0\<s\<2n-2$.
\item[(3)] $\xi T_{[s]}\subseteq T_{[s-1]}$ for all $0\<s\<2n-2$.
\item[(4)] $T^{\xi}=\oplus_{s=0}^{n-1}(T^{\xi}\cap T_{[s]})$ and ${\rm dim}(T^{\xi}\cap T_{[s]})=1$ if $0\<s\<n-1$.
\end{enumerate}
\end{lemma}

\begin{proof}
This follows from a straightforward verification using the action of $x$ and $\xi$ on the basis $\{m_i\otimes m'_j\}$, together with the nilpotent or non-nilpotent behavior of $\mathcal{D}$ and the structure of projective simple modules.
\end{proof}

Let $0\neq x_s\in T^{\xi}\cap T_{[s]}$ for $0\<s\<n-1$, and set $U_s=\langle x_s\rangle$ be the submodule of $T$ generated by $x_s$.

Similarly to Lemma \ref{5.5}, we have the following lemma.

\begin{lemma}{\label{6.9}}
Let $U$ be a simple submodule of $T$. Then  there exists $0\<s\<n-1$ such that $U=U_s$.
\end{lemma}

\begin{lemma}\label{6.10}
Let $0\<s\<n-1$. Then $x_s\in T_{\phi^s\l\l'}$.
\end{lemma}

\begin{proof}
The proof is similar to that of Lemma \ref{5.6}.
\end{proof}

\begin{lemma}\label{6.11}
If $\l\l'\in I'_n$, then $U_s\cong V(n,\phi^s\l\l')$ for all $0\<s\<n-1$.
\end{lemma}

\begin{proof}
Assume that $\l\l'\in I'_n$ and $0\<s\<n-1$. Then $\phi^s\l\l'\in I'_n$.
By Lemma \ref{6.10},  and an argument similar to the proof of Lemma \ref{6.6}(1), we obtain that $U_s\cong V(n,\phi^s\l\l')$.
\end{proof}

In the remainder of this subsection, unless otherwise stated, we assume $\l\l'\notin I'_n$. Let $I_0:=I''_n$. Then $\l\l'\in I_l$ for some $0\<l\<n-1$. In this case, $(\phi^s\l\l')(a\chi^{-1})=\rho^{l-1-2s}$.

\begin{lemma}\label{6.12}
Let $0\<s\<n-1$ and retain the above notations.
\begin{enumerate}
\item[(1)]If $c(l)\<s\<l$, then $1\<n+l-2s\<n$ and
$$U_s\cong V(n+l-2s, \phi^s\l\l').$$
\item[(2)]If $c(n+l)\<s\<n-1$, then $2\<2n+l-2s\<n$ and
$$U_s\cong V(2n+l-2s, \phi^s\l\l').$$
\end{enumerate}
\end{lemma}

\begin{proof} We only consider the case when $n$ and $l$ are both even, since the proofs for other parity combinations are similar. In this case,  $c(l)=\frac{l}{2}$ and $c(n+l)=\frac{n+l}{2}$.

By Lemma \ref{5.1} and \ref{6.10}, we have
\begin{eqnarray}\eqlabel{eq6.2}
   \xi(x^kx_s)=(k)_\rho((\phi^s\l\l')(\chi)-\rho^{1-k}(\phi^s\l\l')(a))x^{k-1}x_s, \ k\>1.
\end{eqnarray}

(1) First, let $s=\frac{l}{2}$. Then $n+l-2s=n$ and $(\phi^s\l\l')(a\chi^{-1})=\rho^{n-1}$. Hence $\phi^s\l\l'\in I''_n$. By Eq.\equref{eq6.2}, $\xi(x^kx_s)=\a_k(n,\phi^s\l\l')x^{k-1}x_s$ for $k\>1$. Induction on $k$ gives  $\xi(x^kx_s)\neq 0$ for all $1\<k\<n-1$.  Hence, by Lemmas \ref{5.2}(1) and \ref{6.10}
$$U_s\cong V(n, \phi^s\l\l').$$
Next, suppose $\frac{l+2}{2}\<s\<l$. Then $2\<n+l-2s\<n-2$ and $\phi^s\l\l'\in I_{n+l-2s}$.
Again by Eq.\equref{eq6.2},  we have
$$\xi(x^kx_s)=\a_k(n+l-2s,\phi^s\l\l')x^{k-1}x_s, \ k\>1.$$

By induction on $k$, it follows that  $\xi(x^kx_s)\neq 0$ for all $1\<k\<n+l-1-2s$.
For $k=n+l-2s$,  we have  $\xi(x^kx_s)=0$. Since $n\<k+s<2n-2$, Lemma \ref{6.8}(1) yields $$x^kx_s\in T_{[k+s]}\oplus T_{[k+s-n]}.$$
 Using Lemma \ref{6.8}(4) and the fact that $x^kx_s\in T^{\xi}$, we deduce $x^kx_s=\a x_{k+s-n}$ for some $\a\in \Bbbk$.

By Lemma \ref{6.10}, we have
$$(g\g)(x^kx_s)=\a (g\g)x_{k+s-n}=\a(\phi^{k+s-n}\l\l')(g\g)x_{k+s-n}=(\phi^{k+s-n}\l\l')(g\g)x^kx_s$$
 for all $g\g\in G\times\G$. On the other hand, using $\xi x_s=0$ and  $2\<k\<n-2$, we have $\xi^{n-1}(x^jx_s)=0$ for all $0\<j\<k$. Then by $gx=\chi^{-1}(g)xg$ and Proposition \ref{2.1}(d, e), one can check that
 $$(g\g)(x^kx_s)=(\phi^{k+s}\l\l')(g\g)x^kx_s$$
  for all $g\g\in G\times\G$. If $x^kx_s\neq0$, then $\phi^{k+s-n}\l\l'=\phi^{k+s}\l\l'$, which forces $\phi^n=1$. Since  ${\rm ord}(\phi)=mn$, this implies $m=1$. Hence,  $\mathcal D$ is of nilpotent type by \cite[Remark 3.4]{SUN C Z}. By Lemma \ref{6.8}(2, 4), we would then have
  $$x^kx_s\in T_{[k+s]}\cap T^{\xi}=0, $$
 a contradiction. Therefore, $x^kx_s=0$, and  it follows from Lemmas \ref{5.2}(2) and \ref{6.10}. that
 $$U_s\cong V(n+l-2s, \phi^s\l\l').$$

(2)  First,  suppose $s=\frac{n+l}{2}$. Then $2n+l-2s=n$ and $(\phi^s\l\l')(a\chi^{-1})=\rho^{n-1}$. Hence, by an argument  similar to case  (1), we have
$$U_s\cong V(n, \phi^s\l\l').$$
Next,  suppose $\frac{n+l+2}{2}\<s\<n-1$. Then $2\<2n+l-2s\<n-2$ and $\phi^s\l\l'\in I_{2n+l-2s}$. Moreover, $n+1\<2n+l-s\<2n-3$. Using a reasoning analogous to case (1), we conclude $$U_s\cong V(2n+l-2s, \phi^s\l\l').$$
\end{proof}

\begin{theorem}\label{6.16}
Let $\l, \l' \in I_n$ and $T=V(n,\l)\otimes V(n,\l')$.
\begin{enumerate}
\item[(1)] If $\l\l'\in I'_n$, then $T\cong \oplus_{s=0}^{n-1} V(n,\phi^s\l\l').$
\item[(2)] If there exists an integer $0\<l\<n-1$ such that $\l\l'\in I_l$, then
$$\begin{array}{rl}
  T\cong& (\oplus_{c(l)\<s\<l}P(n+l-2s,\phi^s\l\l'))\\
   &\oplus(\oplus_{c(n+l)\<s\<n-1}P(2n+l-2s,\phi^s\l\l')).
\end{array}$$
\end{enumerate}
\end{theorem}

\begin{proof}
(1) Assume that $\l\l'\in I'_n$. Then by Lemma \ref{6.11}, $T$ contains a submodule isomorphic to $\oplus_{k=0}^{n-1} V(n,\phi^k\l\l')$.  Comparing their dimensions, we have $T\cong \oplus_{k=0}^{n-1} V(n,\phi^k\l\l')$.

(2) Suppose there is an integer $0\<l\<n-1$ such that $\l\l'\in I_l$. Using the above notations, we only consider the case where both  $n$ and $l$ are  even,  since the proofs for the remaining cases are similar.

Assume that both $n$ and $l$ are  even. Then, by Lemma \ref{6.12}, the module $T$ contains a semisimple submodule
$$(\oplus_{s=\frac{l}{2}}^lU_s)\oplus(\oplus_{s=\frac{n+l}{2}}^{n-1}U_s)$$
 since $l<\frac{n+l}{2}$.  Because  $T$ is injective (and hence projective), it contains a submodule isomorphic to
$$\begin{array}{rl}
&(\oplus_{s=\frac{l}{2}}^{l}I(U_s))\oplus(\oplus_{s=\frac{n+l}{2}}^{n-1}I(U_s))\\
\cong & V(n,\phi^{\frac{l}{2}}\l\l')\oplus V(n,\phi^{\frac{n+l}{2}}\l\l')\\
&\oplus(\oplus_{\frac{l+2}{2}\<s\<l}P(n+l-2s,\phi^s\l\l'))\\
&\oplus(\oplus_{\frac{n+l+2}{2}\<s\<n-1}P(2n+l-2s,\phi^s\l\l'))\\
\cong&(\oplus_{s=\frac{l}{2}}^lP(n+l-2s,\phi^s\l\l'))
\oplus(\oplus_{s=\frac{n+l}{2}}^{n-1}P(2n+l-2s,\phi^s\l\l')).
\end{array}$$
 By comparing dimensions, we conclude that
 $$T\cong(\oplus_{s=\frac{l}{2}}^lP(n+l-2s,\phi^s\l\l'))
\oplus(\oplus_{s=\frac{n+l}{2}}^{n-1}P(2n+l-2s,\phi^s\l\l')).$$
\end{proof}

\subsection{Tensor  products of  simple modules with  projective modules}\label{600.2}
~

In this subsection, we study the tensor products of simple modules with non-simple indecomposable projective modules. We begin by examing the tensor products of non-projective simple modules with non-simple indecomposable projective modules.
\begin{theorem}\label{6.17}
Let $1\<l, l'<n$, $\l\in I_l$ and $\l'\in I_{l'}$, and set $l_1={\rm min}\{l, l'\}$ and $t=l+l'-n-1$.
\begin{enumerate}
\item[(1)] If $t<0$, then
$$\begin{array}{rl}
&V(l, \l)\ot P(l', \l')\\
\cong&(\oplus_{s=0}^{l_1-1}P(l+l'-1-2s,\phi^s\l\l'))\\
&\oplus(\oplus_{c(l+l'-1)\<s\<l-1}\oplus_{k=0}^1 P(n+l+l'-1-2s,\phi^{s-kn}\l\l')).\\
\end{array}$$
\item[(2)] If $t\>0$, then
$$\begin{array}{rl}
&V(l, \l)\ot P(l', \l')\\
\cong&(\oplus_{c(l+l'-1)\<s\<l-1}\oplus_{k=0}^1 P(n+l+l'-1-2s,\phi^{s-kn}\l\l'))\\
&\oplus(\oplus_{s=c(t)}^{t} 2P(l+l'-1-2s,\phi^{s}\l\l'))
\oplus(\oplus_{s=t+1}^{l_1-1}P(l+l'-1-2s,\phi^{s}\l\l')).\\
\end{array}$$
\end{enumerate}
\end{theorem}

\begin{proof}
Let $W= V(l,\l)\ot P(l',\l')$.

(1) Assume $t<0$. First, suppose $l\<l'$. Let $W':=V(l,\l)\ot {\rm soc}P(l',\l')$. Then $W'$ is a submodule of $W$. Since soc$P(l', \l')\cong V(l',\l')$, it follows from Lemma \ref{6.2}(1) that $$W'\cong V(l, \l)\ot V(l',\l')\cong \oplus_{s=0}^{l-1} V( l+l'-1-2s,\phi^s\l\l').$$
Hence $P(W')$ can be embedded into $W$ as a submodule. We have
$$P(W') \cong \oplus_{s=0}^{l-1}P( l+l'-1-2s,\phi^s\l\l').$$
 Since $1\<l+l'-1-2s\<n-1$ for all $0\<s\<l-1$, it follows that dim$(P( l+l'-1-2s,\phi^s\l\l'))=2n$, and hence dim$P(W')=2nl=$dim$(W)$. This implies
$$W\cong\oplus_{s=0}^{l-1}P(l+l'-1-2s,\phi^s\l\l').$$

Next suppose $l'<l$. Applying $V(l,\l)\ot$ to the exact sequence
$$0\rightarrow V(l',\l')\rightarrow P(l',\l')\rightarrow \Omega^{-1}V(l',\l')\rightarrow 0,$$
we obtain  another exact sequence
$$0\rightarrow V(l,\l)\ot V(l',\l') \rightarrow V(l,\l)\ot P(l',\l') \rightarrow V(l,\l)\ot \Omega^{-1}V(l',\l')\rightarrow 0.$$
Note that $l\<n-l'$ and $l+n-l'-n-1=l-l'-1\>0$. By \cite[Lemma 3.11(2)]{SUN C Z}, we have  soc$(\Omega^{-1}V(l',\l'))\cong V(n-l',\s(\l'))\oplus V(n-l',\s^{-1}(\l'))$. Consequently, by Theorem \ref{6.2}, we have
$$ \begin{array}{rl}
&V(l,\l)\ot{\rm soc}(\Omega^{-1}V(l',\l'))\\
\cong&(\oplus_{s=c(l-l'-1)}^{l-l'-1}P(l+n-l'-1-2s,\phi^{s+l'}\l\l'))\\
&\oplus(\oplus_{s=c(l-l'-1)}^{l-l'-1}P(l+n-l'-1-2s,\phi^{s-n+l'}\l\l'))\\
&\oplus(\oplus_{s=l-l'}^{l-1}V(l+n-l'-1-2s,\phi^{s+l'}\l\l'))\\
&\oplus(\oplus _{s=l-l'}^{l-1} V( l+n-l'-1-2s,\phi^{s-n+l'}\l\l')).\\
\end{array}$$
Let $$\begin{array}{rl}
P:=&(\oplus_{s=c(l-l'-1)}^{l-l'-1}P(n+l-l'-1-2s,\phi^{s+l'}\l\l'))\\
&\oplus(\oplus_{s=c(l-l'-1)}^{l-l'-1}P(n+l-l'-1-2s,\phi^{s-n+l'}\l\l')) \\
\end{array}$$
which is both  projective and injective. Thus,  there exists an epimorphism
$$\phi:W:= V(l,\l)\ot P(l',\l')\rightarrow P$$
such that Ker$\phi$ contains a submodule isomorphic to $V(l,\l)\ot V(l',\l')$. Hence $W\cong$Ker$\phi \oplus P$, and so ${\rm soc}(W)\cong{\rm soc(Ker}(\phi))\oplus{\rm soc}(P)$. By Theorem \ref{6.2}(1),
$$V(l,\l)\ot V(l',\l')\cong \oplus_{s=0}^{l'-1}V( l+l'-1-2s,\phi^s\l\l').$$
Let $U={\rm soc}(V(l,\l)\ot V(l',\l'))\oplus {\rm soc}(P)$. Then
$$\begin{array}{rl}
U\cong&(\oplus_{s=0}^{l'-1}V(l+l'-1-2s,\phi^s\l\l'))\\
&\oplus(\oplus_{s=c(l-l'-1)}^{l-l'-1}V(n+l-l'-1-2s,\phi^{s+l'}\l\l'))\\
&\oplus(\oplus_{s=c(l-l'-1)}^{l-l'-1}V(n+l-l'-1-2s,\phi^{s+l'-n}\l\l')).\\
\end{array}$$
Thus, $P(U)$ is isomorphic to a submodule of $W$. A straightforward computation shows that dim$P(U)=2nl=$dim$W$, and hence
$$\begin{array}{rl}
W\cong P(U)\cong&(\oplus_{s=0}^{l'-1}P(l+l'-1-2s,\phi^s\l\l'))\\
&\oplus(\oplus_{s=c(l-l'-1)}^{l-l'-1}P(n+l-l'-1-2s,\phi^{s+l'}\l\l'))\\
&\oplus(\oplus_{s=c(l-l'-1)}^{l-l'-1}P(n+l-l'-1-2s,\phi^{s+l'-n}\l\l'))\\
\cong&(\oplus_{s=0}^{l'-1}P(l+l'-1-2s,\phi^s\l\l'))\\
&\oplus(\oplus_{s=c(l+l'-1)}^{l-1}P(n+l+l'-1-2s,\phi^s\l\l'))\\
&\oplus(\oplus_{s=c(l+l'-1)}^{l-1}P(n+l+l'-1-2s,\phi^{s-n}\l\l')).\\
\end{array}$$

(2) Assume $t\>0$. First,  suppose  $l\<l'$. Then we have an exact sequence
$$0\rightarrow V(l,\l)\ot \Omega V(l',\l')\rightarrow V(l,\l)\ot P(l',\l')\rightarrow V(l,\l)\ot V(l',\l')\rightarrow0.$$
By Theorem \ref{6.2}(2), $P:=\oplus _{s=c(t)}^{t}  P(l+l'-1-2s,\phi^s\l\l')$ is isomorphic to a summand of $V(l,\l)\ot V(l',\l')$. Hence
there exists a module epimorphism
$$\phi:W:= V(l,\l)\ot P(l',\l')\rightarrow P$$
such that Ker$\phi$ contains a submodule isomorphic to $V(l,\l)\ot \Omega V(l',\l')$. Note that
$$V(l,\l)\ot \Omega V(l',\l') \supseteq V(l,\l)\ot\ \text{soc}\  (\Omega V(l',\l'))\cong V(l,\l)\ot V(l',\l').$$
Thus,  by Proposition \ref{6.2}(2), an argument similar to that in (1) show that soc$W$ contains a submodule $U$ isomorphic to
$$\begin{array}{rl}
{\rm soc}(V(l,\l)\ot V(l',\l'))\oplus{\rm soc}(P)
\cong & (\oplus_{s=c(t)}^{t}2V(l+l'-1-2s,\phi^s\l\l'))\\
&\oplus(\oplus_{t+1\<s\<l-1}V(l+l'-1-2s,\phi^s\l\l')).
\end{array}$$
Thus, $P(U)$ is isomorphic to a submodule of $W$. A straightforward computation shows that ${\rm dim}P(U) =2nl ={\rm dim}(W)$. It follows that
$$\begin{array}{rl}
W\cong P(U)\cong&(\oplus_{s=c(t)}^{t}2P(l+l'-1-2s,\phi^s\l\l'))\\
&\oplus(\oplus_{s=t+1}^{l-1}P(l+l'-1-2s,\phi^s\l\l')).\\
\end{array}$$

Next,  suppose  $l'<l$. Then we have two exact sequences:
$$0\rightarrow V(l,\l)\ot V(l',\l')\rightarrow V(l,\l)\ot P(l',\l') \rightarrow V(l,\l)\ot \Omega^{-1} V(l',\l')\rightarrow0,$$
$$\begin{array}{rl}
 0\rightarrow &V(l,\l)\ot V(n-l',\s(\l'))\oplus V(l,\l)\ot V(n-l',\s^{-1}(\l'))\\
&\rightarrow V(l,\l)\ot \Omega^{-1} V(l',\l')\rightarrow V(l,\l)\ot V(l',\l')\rightarrow0.\end{array}$$

Note that $n-l'<l$ and $l+n-l'-n-1=l-l'-1\>0$. By Theorem \ref{6.2}(2), $V(l,\l)\ot V(n-l',\s(\l'))\oplus V(l,\l)\ot V(n-l',\s^{-1}(\l'))$ contains a summand isomorphic to
$$\begin{array}{rl}
   &(\oplus _{s=c(l-l'-1)}^{l-l'-1}P(n+l-l'-1-2s,\phi^{s+l'}\l\l'))\\
   &\oplus (\oplus _{s=c(l-l'-1)}^{l-l'-1}P(n+l-l'-1-2s,\phi^{s+l'-n}\l\l')) \\
   \end{array}$$
and $V(l,\l)\ot V(l',\l')$ contains a summand isomorphic to
$$\oplus _{s=c(t)}^{t}  P(l+l'-1-2s,\phi^s\l\l').$$
It follows  that $V(l,\l)\ot \Omega^{-1} V(l',\l')$ contains a projective summand $\overline{P}$ isomorphic to
$$\begin{array}{rl}
  &(\oplus_{s=c(l-l'-1)}^{l-l'-1}P(n+l-l'-1-2s,\phi^{s+l'}\l\l'))\\
  &\oplus(\oplus_{s=c(l-l'-1)}^{l-l'-1}P(n+l-l'-1-2s,\phi^{s+l'-n}\l\l')) \\
   &\oplus(\oplus_{s=c(t)}^{t}P(l+l'-1-2s,\phi^{s}\l\l')).\end{array}$$
Then, an argument similar to the proof of (1) shows that soc$(V(l,\l)\ot P(l',\l'))$ contains a submodule $U$ isomorphic to
${\rm soc}(V(l,\l)\ot V(l',\l'))\oplus {\rm soc}(\overline{P})$. By Corollary \ref{5.8}(2) or Theorem \ref{6.2}(2), we have
$$\begin{array}{rl}
U\cong&\oplus_{s=c(t)}^{l'-1} V(l+l'-1-2s,\phi^s\l\l')\oplus{\rm soc}\overline{P}\\
\cong&(\oplus_{s=c(l-l'-1)}^{l-l'-1}V(n+l-l'-1-2s,\phi^{s+l'}\l\l'))\\
&\oplus(\oplus_{s=c(l-l'-1)}^{l-l'-1}V(n+l-l'-1-2s,\phi^{s+l'-n}\l\l'))\\
&\oplus(\oplus_{s=c(t)}^{t}2V(l+l'-1-2s,\phi^{s}\l\l'))\\
&\oplus(\oplus_{s=t+1}^{l'-1} V(l+l'-1-2s,\phi^{s}\l\l')).\\
\end{array}$$
Then one can check that ${\rm dim}P(U)={\rm dim}W$, and so
$$\begin{array}{rl}
W\cong&(\oplus_{s=c(l-l'-1)}^{l-l'-1}P(n+l-l'-1-2s,\phi^{s+l'}\l\l'))\\
&\oplus(\oplus_{s=c(l-l'-1)}^{l-l'-1}P(n+l-l'-1-2s,\phi^{s+l'-n}\l\l'))\\
&\oplus(\oplus_{s=c(t)}^{t}2P(l+l'-1-2s,\phi^{s}\l\l'))\\
&\oplus(\oplus_{s=t+1}^{l'-1}P(l+l'-1-2s,\phi^{s}\l\l'))\\
\cong&(\oplus_{s=c(l+l'-1)}^{l-1}P(n+l+l'-1-2s,\phi^{s}\l\l'))\\
&\oplus(\oplus_{s=c(l+l'-1)}^{l-1}P(n+l+l'-1-2s,\phi^{s-n}\l\l'))\\
&\oplus(\oplus_{s=c(t)}^{t}2P(l+l'-1-2s,\phi^{s}\l\l'))\\
&\oplus(\oplus_{s=t+1}^{l'-1}P(l+l'-1-2s,\phi^{s}\l\l')).\\
\end{array}$$
\end{proof}

Next, we consider the tensor products of simple projective modules with non-simple indecomposable projective modules.
\begin{theorem}\label{6.18}
Let $1\<l'\<n-1$, $\l \in I_n$ and $\l'\in I_{l'}$.
\begin{enumerate}
\item[(1)] If $\l\in I'_n$, then
$$V(n,\l)\ot P(l',\l')
\cong(\oplus_{s=0}^{l'-1}2V(n,\phi^s\l\l'))
\oplus(\oplus_{s=l'}^{n-1}\oplus_{k=0}^1 V(n,\phi^{s-kn}\l\l')).$$
\item[(2)] If $\l\in I''_n$, then
$$\begin{array}{rl}
V(n,\l)\otimes P(l',\l')
\cong&(\oplus_{s=c(l'-1)}^{l'-1}2P(n+l'-1-2s,\phi^s\l\l'))\\
&\oplus(\oplus_{s=c(n+l'-1)}^{n-1}\oplus_{k=0}^1 P(2n+l'-1-2s,\phi^{s-kn}\l\l')).
\end{array}$$
\end{enumerate}
\end{theorem}

\begin{proof}
 (1) Assume $\l \in I'_n$. Since $V(n,\l)$ is projective,  by \cite[Lemma 3.11]{SUN C Z}, we have
$$\begin{array}{rl}
&V(n,\l)\ot P(l',\l')\\
\cong &2(V(n,\l)\ot V(l',\l'))\oplus(V(n,\l)\ot V(n-l',\s(\l')))\\
&\oplus(V(n,\l)\ot V(n-l,\s^{-1}(\l'))).\end{array}$$
 Then,  by Theorem \ref{6.7}(1), we have
\begin{eqnarray*}
&V(n,\l)\ot P(l',\l')&
\cong(\oplus_{s=0}^{l'-1}2V(n,\phi^s\l\l')) \oplus
   % (\oplus_{s=0}^{n-l'-1}V(n,\phi^{s+l'}\l\l') \oplus
    (\oplus_{s=0}^{n-l'-1}\oplus_{k=0}^1V(n,\phi^{s+l'-kn}\l\l'))\\
&&\cong(\oplus_{s=0}^{l'-1}2V(n,\phi^s\l\l')) \oplus
   % (\oplus_{s=0}^{n-l'-1}V(n,\phi^{s+l'}\l\l') \oplus
    (\oplus_{s=l'}^{n-1}\oplus_{k=0}^1V(n,\phi^{s-kn}\l\l')).\\
\end{eqnarray*}
(2) The proof is analogous to that of $(1)$, and is therefore omitted.
\end{proof}

\subsection{The tensor product of two projective modules}
~

In this subsection, we study the tensor products of two non-simple indecomposable projective modules.

\begin{theorem}\label{6.19}
Let $1\<l,l'\<n-1$, $\l\in I_l$ and $\l'\in I_{l'}$, and let $l_1$={\rm min}$\{l,l'\}$, $l_2$={\rm max}$\{l,l'\}$ and $t=l+l'-(n+1)$.
\begin{enumerate}
\item[(1)] If $t<0$, then
\begin{eqnarray*}
&&P(l,\l)\ot P(l',\l')\\
&\cong&(\oplus_{s=0}^{l_1-1}2P(l+l'-1-2s,\phi^s\l\l'))\\
&&\oplus(\oplus_{c(l+l'-1)\<s\<l_2-1}\oplus_{k=0}^1 2P(n+l+l'-1-2s,\phi^{s-kn}\l\l'))\\
&&\oplus(\oplus_{s=l_2}^{l+l'-1}\oplus_{k=0}^1 P(n+l+l'-1-2s,\phi^{s-kn}\l\l'))\\
&&\oplus(\oplus_{c(n+l+l'-1)\<s\<n-1}\oplus_{k,p=0}^1 P(2n+l'+l-1-2s,\phi^{s-kn-pn}\l\l')).\\
\end{eqnarray*}
\item[(2)] If $t\>0$, then
\begin{eqnarray*}
& & P(l,\l)\ot P(l',\l')\\
&\cong&(\oplus_{s=c(t)}^{t}4P(l+l'-1-2s,\phi^{s}\l\l'))
\oplus(\oplus_{s=t+1}^{l_1-1}2P(l+l'-1-2s,\phi^s\l\l'))\\
&&\oplus(\oplus_{c(l+l'-1)\<s\<l_2-1}\oplus_{k=0}^12P(n+l+l'-1-2s,\phi^{s-kn}\l\l'))\\
&&\oplus(\oplus_{s=l_2}^{n-1}\oplus_{k=0}^1P(n+l+l'-1-2s,\phi^{s-kn}\l\l')).\\
\end{eqnarray*}
\end{enumerate}
\end{theorem}

\begin{proof}
Since $P(l,\l)\ot P(l',\l')\cong P(l',\l')\ot P(l,\l)$, we may assume that $l\>l'$. By \cite[Lemma 3.11]{SUN C Z}, we have
$$P(l,\l)\ot P(l',\l')\cong 2(V(l,\l)\ot P(l',\l'))\oplus (\oplus_{k=0}^1V(n-l,(\tau^{-k}\s)(\l))\ot P(l',\l')).$$
Moreover, $n-l+l'-(n+1)=l'-l-1<0$.

(1) Assume $t<0$. Then $n-l\>l'$. By Proposition \ref{6.17}(1),
$$\begin{array}{rl}
&P(l,\l)\ot P(l',\l')\\
\cong &(\oplus_{s=0}^{l'-1}2P(l+l'-1-2s,\phi^s\l\l'))\\
&\oplus(\oplus_{c(l+l'-1)\<s\<l-1}\oplus_{k=0}^1 2P(n+l+l'-1-2s,\phi^{s-kn}\l\l'))\\
&\oplus(\oplus_{k=0}^1\oplus_{s=0}^{l'-1}P(n-l+l'-1-2s,\phi^{l-kn+s}\l\l'))\\
&\oplus(\oplus_{c(n-l+l'-1)\<s\<n-l-1}\oplus_{k,p=0}^1P(2n-l+l'-1-2s,\phi^{l-kn-pn+s}\l\l'))\\
\cong&(\oplus_{s=0}^{l'-1}2P(l+l'-1-2s,\phi^s\l\l'))\\
&\oplus(\oplus_{c(l+l'-1)\<s\<l-1}\oplus_{k=0}^12P(n+l+l'-1-2s,\phi^{s-kn}\l\l'))\\
&\oplus(\oplus_{s=l}^{l+l'-1}\oplus_{k=0}^1P(n+l+l'-1-2s,\phi^{s-kn}\l\l'))\\
&\oplus(\oplus_{c(n+l+l'-1)\<s\<n-1}\oplus_{k,p=0}^1P(2n+l+l'-1-2s,\phi^{s-kn-pn}\l\l')).
\end{array}$$

(2) Assume $t\>0$. Then $n-l<l'$, and hence $c(n-l+l'-1)>n-l-1$. Thus, by Proposition \ref{6.17}, we have
$$\begin{array}{rl}
&P(l,\l)\ot P(l',\l')\\
\cong &(\oplus_{s=c(t)}^t 4P(l+l'-1-2s,\phi^s\l\l'))\oplus(\oplus_{s=t+1}^{l'-1}2P(l+l'-1-2s,\phi^s\l\l'))\\
&\oplus(\oplus_{c(l+l'-1)\<s\<l-1}\oplus_{k=0}^1 2P(n+l+l'-1-2s,\phi^{s-kn}\l\l'))\\
&\oplus(\oplus_{k=0}^1 \oplus_{s=0}^{n-l-1}P(n-l+l'-1-2s,\phi^{l-kn+s}\l\l'))\\
\cong &(\oplus_{s=c(t)}^t4P(l+l'-1-2s,\phi^{s}\l\l')) \oplus(\oplus_{s=t+1}^{l'-1}2P(l+l'-1-2s,\phi^{s}\l\l'))\\
&\oplus(\oplus_{c(l+l'-1)\<s\<l-1}\oplus_{k=0}^12P(n+l+l'-1-2s,\phi^{s-kn}\l\l'))\\
&\oplus(\oplus_{s=l}^{n-1}\oplus_{k=0}^1P(n+l+l'-1-2s,\phi^{s-kn}\l\l')).
\end{array}$$
\end{proof}

\section{The projecitve class rings $r_p(D(H_{\mathcal{D}}))$}\label{s700}

In this section, we investigate the Grothendieck ring $G_0(D(H_{\mathcal{D}}))$ and projective class ring $r_p(D(H_{\mathcal{D}}))$.
\subsection{The Grothendieck ring $G_0(D(H_{\mathcal{D}}))$}
~

In this subsection, we compute the  Grothendieck ring $G_0(D(H_{\mathcal{D}}))$. Note that $\{[V(l,\l)]|1\<l\<n,\l\in I_l\}$ is a $\mathbb{Z}$-basis of $G_0(D(H_{\mathcal{D}}))$.
\begin{proposition}\label{7.1}
Assume that $1\<l\<l'<n$, $\l\in I_{l}$ and $\l'\in I_{l'}$.
Let $t=l+l'-n-1$. In $G_0(D(H_{\mathcal{D}}))$, we have
\begin{enumerate}
\item[(1)] If  $t<0$, then $[V(l,\l)] [V(l',\l')]=\sum _{s=0}^{l-1}[V( l+l'-1-2s,\phi^s\l\l')]$.
\item[(2)] If $t\>0$, then

\begin{equation*}
\begin{array}{ll}
\vspace{0.2 cm}
[V(l,\l)][V(l',\l')]=\left\{
\begin{array}{ll}
H, & {\rm if}\ t \ {\rm is \ odd},\\
H+[V(n, \phi^{\frac{t}{2}}\l\l')], &{\rm if}\ t\ {\rm is \ even},\\
\end{array}\right.\\
\end{array}
\end{equation*}
where

 $$\begin{array}{rl}
H= &\sum_{t+1\<s\<l-1}[V(l+l'-1-2s,\phi^s\l\l')]\\
&+\sum_{c(t+1)\<s\<t}2[V(l+l'-1-2s,\phi^s\l\l')]\\
&+\sum_{c(t+1)\<s\<t}\sum_{i=0}^1[V(n-l-l'+1+2s,\phi^{l+l'-1-s-in}\l\l')].\end{array}$$

\end{enumerate}
\end{proposition}
\begin{proof}
The result follows directly  from Theorem \ref{6.2} together with \cite[Lemma 3.11]{SUN C Z}.
\end{proof}

\begin{proposition}\label{7.2}
Let $1\<l<n$, $\l\in I_l$ and $\l'\in I_n$. In $G_0(D(H_{\mathcal{D}}))$, we have
\begin{enumerate}
\item[(1)] If $\l'\in I'_n$, then $[V(l,\l)][V(n,\l')]=\sum_{s=0}^{l-1} [V(n,\phi^s\l\l')]$.
\item[(2)] If $\l'\in I''_n$, then

\begin{equation*}
\begin{array}{ll}
\vspace{0.2 cm}
[V(l,\l)][V(n,\l')]=\left\{
\begin{array}{ll}
W+[V(n,\phi^{\frac{l-1}{2}}\l\l')], & {\rm if} \ l \ {\rm is \ odd},\\
W, & {\rm if} \ l \ {\rm is \ even},\\
\end{array}\right.\\
\end{array}
\end{equation*}
where
 $$\begin{array}{rl}
 W=&\sum _{s=c(l)}^{l-1} 2[V(n+l-1-2s,\phi^{s}\l\l')]\\
 &+\sum _{s=c(l)}^{l-1} \sum_{i=0}^1[V(2s+1-l,\phi^{in+l-1-s}\l\l')].\\
    \end{array}$$
    \end{enumerate}
\end{proposition}
\begin{proof}
It follows from Theorem \ref{6.7} and \cite[Lemma 3.11]{SUN C Z}.
\end{proof}

\begin{proposition}\label{7.3}
Let $\l, \l' \in I_n$. In $G_0(D(H_{\mathcal{D}}))$, we have
\begin{enumerate}
\item[(1)] If $\l\l'\in I'_n$, then $[V(n,\l)][V(n,\l')]= \sum_{s=0}^{n-1} [V(n,\phi^s\l\l')].$
\item[(2)] If there is an integer $0\<l\<n-1$ such that $\l\l'\in I_l$, then
\begin{equation*}
\begin{array}{ll}
\vspace{0.2 cm}
[V(n,\l)][V(n,\l')]=\left\{
\begin{array}{ll}
U+[V(n, \phi^{\frac{n+l}{2}}\l\l')],& {\rm if} \ l \ {\rm and}\  n \ {\rm are \  odd},\\
U, & {\rm if}\ l \ {\rm is \ odd} \ {\rm and} \ n \ {\rm is \ even},\\
U+[V(n, \phi^{\frac{l}{2}}\l\l')],& {\rm if } \ l \ {\rm is \ even} \ {\rm and}  \ n \ {\rm is \ odd},\\
U+\sum_{i=0}^1[V(n, \phi^{\frac{in+l}{2}}\l\l')], & {\rm if} \ l \ {\rm and} \ n \ {\rm are \ even},\\
\end{array}\right.\\
\end{array}
\end{equation*}
where $I_0=I''_n$ and
$$\begin{array}{rl}
U=& \sum_{c(l+1)\<s\<l}2[V(n+l-2s,\phi^s\l\l')]\\
   &+\sum_{c(l+1)\<s\<l}\sum_{i=0}^1[V(2s-l,\phi^{in+l-s}\l\l')]\\
  & +\sum_{c(n+l+1)\<s\<n-1}2[V(2n+l-2s,\phi^s\l\l')]\\
   &+\sum_{c(n+l+1)\<s\<n-1}\sum_{i=1}^2[V(2s-l-n,\phi^{in+l-s}\l\l')].
\end{array}$$
\end{enumerate}
\end{proposition}
\begin{proof}
It follows from Theorem \ref{6.16} and \cite[Lemma 3.11]{SUN C Z}.
\end{proof}

\begin{lemma}\label{7.4}
Let $2\<s\<n-1$. Then
$$V(2,\chi)^{\ot s}\cong \oplus_{i=0}^{[\frac{s}{2}]}\frac{s-2i+1}{s-i+1}\binom{s}{i}V(s+1-2i,\phi^i\chi^s).$$
\end{lemma}
\begin{proof}
The proof  is similar to that of \cite[Lemma 5.3]{ChHasLinSun}.
\end{proof}

Let $x_{\l}=[V(1,\l)],\l\in K$, $y=[V(2,\chi)]$ and $z_{\l'}=[V(n,\l')]$, $\l'\in I'_n$.

\begin{lemma}\label{7.5}
Let $1\<l'<n$ with $\varphi\in I_{l'}$ and let $\psi\in I''_n$. The following equations hold in $G_0(D(H_D))$:
\begin{enumerate}
\item[(1)]
$[V(l',\varphi)]=[V(1,\varphi\chi^{1-l'})][V(l',\chi^{l'-1})]$,
\item[(2)]
$[V(n,\psi)]=[V(1,\psi\chi^{1-n})][V(n,\chi^{n-1})]$,
\item[(3)]
$[V(l,\chi^{l-1})]=\sum_{i=0}^{[\frac{l-1}{2}]}(-1)^i\binom{l-1-i}{i}
x^i_{\phi\chi^2}y^{l-1-2i}$, $1\<l\<n$.
\end{enumerate}
\end{lemma}
\begin{proof}
(1) and (2) follow directly from Lemma \ref{5.3} and the equalities $\sharp I_l=\sharp I''_n=\sharp K$.

(3) We prove the statement by induction on $l$. For $l=1$, $l=2$ or $l=3$, the equation holds trivially.  Assume now  that $2\<l\<n-1$ and the statement holds for all smaller values of $l$.  By Theorem \ref{6.2}, we have
$$V(2,\chi)\ot V(l,\chi^{l-1})\cong V(l+1,\chi^l)\oplus V(l-1,\phi\chi^l).$$
By the induction hypothesis, we have

$$\begin{array}{rl}
[V(l+1,\chi^l)]=&y[V(l,\chi^{l-1})]-[V(l-1,\chi^{l-2})]x_{\phi\chi^2}\\
=&\sum_{i=0}^{[\frac{l-1}{2}]}(-1)^i\binom{l-1-i}{i}x^i_{\phi\chi^2}y^{l-2i}-\\
&\sum_{i=0}^{[\frac{l-2}{2}]}(-1)^i\binom{l-2-i}{i}x^{i+1}_{\phi\chi^2}y^{l-2-2i}\\  =&\sum_{i=0}^{[\frac{l}{2}]}(-1)^i\binom{l-i}{i})x^i_{\phi\chi^2}y^{l-2i}.\\
\end{array}$$
This completes the proof.
\end{proof}

\begin{corollary}\label{7.6}
The commutative ring $G_0(D(H_{\mathcal{D}}))$ is generated by $\{x_{\l},y, z_{\l'}|\l\in K, \l'\in I'_n\}$.
\end{corollary}
\begin{proof}
Let $R$ be the subring of $G_0(D(H_{\mathcal{D}}))$ generated by  $\{x_{\l},y, z_{\l'}|\l\in K, \l'\in I'_n\}$. Clearly,  $R\subseteq G_0(D(H_{\mathcal{D}}))$.  By Lemma \ref{7.5}, we have $[V(l,\varphi)]\in R$ and $[V(n,\psi)]\in R$ for $1\<l<n$, $\varphi\in I_l$ and $\psi\in I''_n$. Hence,  $G_0(D(H_{\mathcal{D}})\subseteq R$. Therefore,  $G_0(D(H_{\mathcal{D}})= R$.
\end{proof}

\begin{lemma}\label{7.7}
In $G_0(D(H_{\mathcal{D}}))$, we have
$$\sum_{i=0}^{[\frac{n}{2}]}(-1)^i\frac{n}{n-i}\binom{n-i}{i}x^i_{\phi\chi^2}y^{n-2i}-x_{\phi^n\chi^n}-x_{\chi^n}=0.$$
\end{lemma}
\begin{proof}
We only consider the case when $n$ is odd, the proof for even $n$ is similar.

By Lemma \ref{7.5}(3), we have
$$[V(n,\chi^{n-1})]=\sum_{i=0}^{\frac{n-1}{2}}(-1)^i\binom{n-1-i}{i}
x^i_{\phi\chi^2}y^{n-1-2i}.$$
Hence,
$$y[V(n,\chi^{n-1})]=\sum_{i=0}^{\frac{n-1}{2}}(-1)^i\binom{n-1-i}{i}
x^i_{\phi\chi^2}y^{n-2i}.$$
 On the other hand, by Proposition \ref{7.2} and Lemma \ref{7.5}(1) we also have
$$y[V(n,\chi^{n-1})]=2[V(n-1,\chi^{n-2})]x_{\phi\chi^2}+x_{\phi^n\chi^n}+x_{\chi^n}.$$
Combing these two expressions yields
$$\sum_{i=0}^{\frac{n-1}{2}}(-1)^i\binom{n-1-i}{i}
x^i_{\phi\chi^2}y^{n-2i}=2\sum_{i=0}^{\frac{n-3}{2}}(-1)^i\binom{n-2-i}{i}x^{i+1}_{\phi\chi^2}y^{n-2-2i}+x_{\phi^n\chi^n}+x_{\chi^n}.$$
A direct verification shows that this equation is  equivalent to
$$\sum_{i=0}^{\frac{n-1}{2}}(-1)^i\frac{n}{n-i}\binom{n-i}{i}x^i_{\phi\chi^2}y^{n-2i}=x_{\phi^n\chi^n}+x_{\chi^n}.$$
This completes the proof.
\end{proof}
For any $1\<l\<n$, let $f_l(x_{\phi\chi^2},y)=\sum_{i=0}^{[\frac{l-1}{2}]}(-1)^i\binom{l-1-i}{i}x^i_{\phi\chi^2}y^{l-1-2i}$ in $G_0(D(H_{\mathcal{D}}))$.
\begin{proposition}\label{7.8}
Let $\l,\mu\in K$ and $\l',\l''\in I'_n$.
\begin{enumerate}
\item[(1)] $x_{\l}x_{\mu}=x_{\l\mu}.$
\item[(2)] $x_{\l}z_{\l'}=z_{\l\l'}.$
\item[(3)] $yz_{\l'}=z_{\chi\l'}+z_{\phi\chi\l'}.$
\item[(4)]  If $\l'\l''\in I'_n$, then $z_{\l'}z_{\l''}=\sum_{j=0}^{n-1}z_{\phi^j\l'\l''}.$
\item[(5)]  If $\l'\l''\in I_l$ for some $0\<l\<n-1$, then
    \end{enumerate}

    \begin{equation*}
\begin{array}{ll}
\vspace{0.2 cm}
z_{\l'}z_{\l''}=\left\{
\begin{array}{ll}
U_{l,\l',\l''}+f_n(x_{\phi\chi^2},y)x_{\phi^{\frac{n+l}{2}}\l'\l''\chi^{1-n}},& {\rm if }\ l \ {\rm and}\  n \ {\rm are\  odd},\\
U_{l,\l',\l''}, & {\rm if}\ l\ {\rm is \ odd}\ {\rm and} \ n \ {\rm is \ even},\\
U_{l,\l',\l''}+f_n(x_{\phi\chi^2},y)x_{\phi^{\frac{l}{2}}\l'\l''\chi^{1-n}},& {\rm if} \ l \ {\rm is \ even} \ {\rm and}  \ n \ {\rm is \ odd},\\
U_{l,\l',\l''}+\sum_{i=0}^1f_n(x_{\phi\chi^2},y)x_{\phi^{\frac{in+l}{2}}\l'\l''\chi^{1-n}}, & {\rm if}\ l \ {\rm and} \ n \ {\rm are \ even},\\
\end{array}\right.\\
\end{array}
\end{equation*}

where
$$\begin{array}{rl}
U_{l,\l',\l''}=& \sum_{c(l+1)\<s\<l}2f_{n+l-2s}(x_{\phi\chi^2},y)x_{\phi^s\l'\l''\chi^{2s+1-l-n}}\\
   &+\sum_{c(l+1)\<s\<l}\sum_{i=0}^1f_{2s-l}(x_{\phi\chi^2},y)x_{\phi^{in+l-s}\l'\l''\chi^{1+l-2s}}\\
  & +\sum_{c(n+l+1)\<s\<n-1}2f_{2n+l-2s}(x_{\phi\chi^2},y)x_{\phi^s\l'\l''\chi^{1+2s-l-2n}}\\
   &+\sum_{c(n+l+1)\<s\<n-1}\sum_{i=1}^2f_{2s-l-n}(x_{\phi\chi^2},y)x_{\phi^{in+l-s}\l'\l''\chi^{1+l+n-2s}}.
\end{array}$$
\end{proposition}
\begin{proof}
(1) follows from Proposition \ref{7.1}. (2) and (3) follow from Proposition \ref{7.2}(1). (4) follows from Proposition \ref{7.3}(1). (5) follows from Proposition \ref{7.3}(2) and Lemma \ref{7.5}.
\end{proof}

\begin{corollary}\label{7.9} The set
$$\{x_{\l}y^i,z_{\l'}|0\<i\<n-1,\l\in K,\l'\in I'_n\}$$ is a $\mathbb{Z}$-basis of $G_0(D(H_{\mathcal{D}}))$.
\end{corollary}
\begin{proof}
By Lemma \ref{7.7}, we have
$$y^n=x_{\phi^n\chi^n}+x_{\chi^n}-\sum_{i=1}^{[\frac{n}{2}]}(-1)^i
\frac{n}{n-i}\binom{n-i}{i}x^i_{\phi\chi^2}y^{n-2i}.$$
 Thus, it follows from Proposition \ref{7.8} that $G_0(D(H_{\mathcal{D}}))$ is generated, as
 $\mathbb{Z}$-module,  by
 $$\{x_{\l}y^i,z_{\l'}|0\<i\<n-1,\l\in K,\l'\in I'_n\}.$$
Since the rank of $\mathbb{Z}$-module $G_0(D(H_{\mathcal{D}}))$ equals the cardinality of this set,   $\{x_{\l}y^i,z_{\l'}|0\<i\<n-1,\l\in K,\l'\in I'_n\}$ forms a $\mathbb{Z}$-basis of $G_0(D(H_{\mathcal{D}}))$.
\end{proof}

Now we proceed to describe the structure of the Grothendieck ring $G_0(D(H_{\mathcal{D}}))$ separately for the cases when $n$ is even and when $n$ is odd.
 Let $X=\{x_{\l},y, z_{\l'}|\l\in K, \l'\in I'_n\}$ and denote by $\mathbb{Z}[X]$  the corresponding polynomial ring.

\textbf{Case 1: $n$ is even}

In this case, define a subset $A_0\subset \mathbb{Z}[X]$ by
$$A_0=\left\{\left.\begin{array}{l}
x_{\l}x_{\mu}-x_{\l\mu};\  x_{\l}z_{\l'}-z_{\l\l'};\\
yz_{\l'}-z_{\chi\l'}-z_{\phi\chi\l'};\\
z_{\l_1}z_{\l_2}-\sum_{j=0}^{n-1}z_{\phi^j\l_1\l_2};\\
z_{\l_3}z_{\l_4}-U_{l,\l_3,\l_4};\\
z_{\l_5}z_{\l_6}-U_{k,\l_5,\l_6}-\\
\sum_{i=0}^1f_n(x_{\phi\chi^2},y)x_{\phi^{\frac{in+k}{2}}\l_5\l_6\chi^{1-n}};\\
\sum_{i=0}^{\frac{n}{2}}(-1)^i\frac{n}{n-i}\binom{n-i}{i}x^i_{\phi\chi^2}y^{n-2i}\\
-x_{\phi^n\chi^n}-x_{\chi^n}.\\
 \end{array}
\right|\begin{array}{l}
\l,\mu\in K;\\
\l',\l_1,\l_2\in I'_n
{\rm with} \ \l_1\l_2\in I'_n;\\
\l_3,\l_4\in I'_n
{\rm with} \ \l_3\l_4\in I_l\\
{\rm for \ some \ odd} \ l , 0\<l\<n-1;\\
\l_5,\l_6\in I'_n
{\rm with} \ \l_5\l_6\in I_k\\
{\rm for \ some \ even} \ k, 0\<k\<n-1.\\
\end{array}\right\},$$
where $U_{l,\l_3,\l_4}$ and $U_{k,\l_5,\l_6}$ are given as in Proposition \ref{7.8}.

\begin{theorem}\label{7.10}
Assume that $n$ is even. Then the Grothendieck ring $G_0(D(H_{\mathcal{D}}))$ is isomorphic to the quotient ring  $\mathbb{Z}[X]/( A_0)$, where $(A_0)$ denotes the ideal of $\mathbb{Z}[X]$ generated by $A_0$.
\end{theorem}
\begin{proof}
By Corollary \ref{7.6}, there is a ring epimorphism $F:\mathbb{Z}[X]\rightarrow G_0(D(H_{\mathcal{D}}))$ defined by
$$F(x_{\l})=[V(1,\l)],\ F(y)=[V(2,\chi)],\  F(z_{\l'})=[V(n,\l')]$$
 for $\l\in K$ and $\l'\in I'_n$. By Lemma \ref{7.7} and Proposition \ref{7.8},  we have $F(A_0)=0$. Hence,  $F$ induces a ring epimorphism
 $$\overline{F}:\  \mathbb{Z}[X]/(A_0)\rightarrow G_0(D(H_{\mathcal{D}}))$$
 such that $F=\overline{F}\pi$, where $\pi:\mathbb{Z}[X]\rightarrow \mathbb{Z}[X]/(A_0)$ is the canonical projection.  By the construction of $A_0$, the quotient $\mathbb{Z}[X]/(A_0)$ is generated, as a $\mathbb{Z}$-module, by
 $$\{\pi(x_{\l})\pi (y)^i,\pi (z_{\l'})|0\<i\<n-1,\l\in K,\l'\in I'_n\}.$$
It follows from  Corollary \ref{7.9}  that $\overline{F}$ is a $\mathbb{Z}$-module isomorphism.  Consequently, $\overline{F}$ is a ring isomorphism.
\end{proof}

\textbf{Case 2: $n$ is odd}

In this case, define a  subset $A_1\subset \mathbb{Z}[X]$ by

$$A_1=\left\{\left.\begin{array}{l}
x_{\l}x_{\mu}-x_{\l\mu};\  x_{\l}z_{\l'}-z_{\l\l'};\\
yz_{\l'}-z_{\chi\l'}-z_{\phi\chi\l'};\\
z_{\l_1}z_{\l_2}-\sum_{j=0}^{n-1}z_{\phi^j\l_1\l_2};\\
z_{\l_3}z_{\l_4}-U_{l,\l_3,\l_4}-\\
f_n(x_{\phi\chi^2},y)x_{\phi^{\frac{n+l}{2}}\l_3\l_4\chi^{1-n}};\\
z_{\l_5}z_{\l_6}-U_{k,\l_5,\l_6}-\\
f_n(x_{\phi\chi^2},y)x_{\phi^{\frac{k}{2}}\l_5\l_6\chi^{1-n}};\\
\sum_{i=0}^{\frac{n-1}{2}}(-1)^i\frac{n}{n-i}\binom{n-i}{i}x^i_{\phi\chi^2}y^{n-2i}\\
-x_{\phi^n\chi^n}-x_{\chi^n}.\\
 \end{array}
\right|\begin{array}{l}
\l,\mu\in K;\\
\l',\l_1,\l_2\in I'_n
{\rm with} \ \l_1\l_2\in I'_n;\\
\l_3,\l_4\in I'_n
{\rm with} \ \l_3\l_4\in I_l\\
{\rm for \ some \ odd} \ l , 0\<l\<n-1;\\
\l_5,\l_6\in I'_n
{\rm with} \ \l_5\l_6\in I_k\\
{\rm for \ some \ even} \ k, 0\<k\<n-1.\\
\end{array}\right\},$$
where $U_{l,\l_3,\l_4}$ and $U_{k,\l_5,\l_6}$ are given as in Proposition 6.8.

\begin{theorem}\label{nodd} Assume that $n$ is odd. Then the Grothendieck ring $G_0(D(H_{\mathcal{D}}))$ is isomorphic to the quotient ring  $\mathbb{Z}[X]/( A_1)$, where  $(A_1)$ is the ideal  of $\mathbb{Z}[X]$ generated by $A_1$.
\end{theorem}
\begin{proof}
The proof  is similar to that of Theorem \ref{7.10}.
\end{proof}

\subsection{The projective class ring $r_p(D(H_{\mathcal{D}}))$}
~

In this subsection, we investigate the  projective class ring $r_p(D(H_{\mathcal{D}}))$. By Section \ref{s500}, the subcategory of  $D(H_{\mathcal{D}})$-mod consisting of semisimple and projective modules forms a monoidal subcategory. Moreover,
$$\{[V(l,\l)],[P(l',\l')]|1\<l\<n,1\<l'\<n-1, \l\in I_l,\l'\in I_{l'}\}$$
forms  a $\mathbb{Z}$-basis of $r_p(D(H_{\mathcal{D}}))$.

For convenience, set
$$x_{\l}=[V(1,\l)], \ y=[{V(2,\chi)}],\ z_{\l'}=[V(n,\l')]$$
 in $r_p(D(H_{\mathcal{D}}))$, where $\l\in K$ and $\l'\in I'_n$. Note that $x_{\l}$ is invertible for any $\l \in K$, with inverse $x^{-1}_{\l}=x_{\l^{-1}}$, by Lemma \ref{5.3}.
\begin{proposition}\label{7.11}
The following equations hold in $r_p(D(H_{\mathcal{D}}))$:
\begin{enumerate}
\item[(1)]
$[V(l,\l)]=[V(1,\l\chi^{1-l})][V(l,\chi^{l-1})]$ for all $1\<l<n$ and $\l\in I_{l}$;
\item[(2)]
$[V(n,\l)]=[V(1,\l\chi^{1-n})][V(n,\chi^{n-1})]$ for all $\l\in I''_n$;
\item[(3)]
$[V(l+1,\chi^{l})]=y^l-\sum_{i=1}^{[\frac{l}{2}]}\frac{l+1-2i}{l+1-i}\binom{l}{i}
x^i_{\phi\chi^2}[V(l+1-2i,\chi^{l-2i})]$ for all $1\<l<n$;
\item[(4)]
$[P(l,\l)]=[V(1,\l\chi^{1-l})][P(l,\chi^{l-1})]$ for all $1\<l<n$ and $\l\in I_l$;
\item[(5)]
$y[V(n,\chi^{n-1})]=x_{\phi\chi^2}[P(n-1,\chi^{n-2})]$;
\item[(6)]
$y[P(1,\varepsilon)]=[P(2,\chi)]+[V(n,\chi^{n-1})](x_{\phi\chi^{2-n}}
+x_{\phi^{1-n}\chi^{2-n}})$;
\item[(7)]
$y[P(n-1,\chi^{n-2})]=2[V(n,\chi^{n-l})]+x_{\phi\chi^2}[P(n-2,\chi^{n-3})]$;
\item[(8)]
$y[P(l,\chi^{l-1})]=[P(l+1,\chi^{l})]+x_{\phi\chi^2}[P(l-1,\chi^{l-2})]$ for all $2\<l\<n-2$.
\end{enumerate}
\end{proposition}
\begin{proof}
Statements (1) and (2) follow directly  from Lemma \ref{5.3}. Statement (3) follows from (1) together with Lemma \ref{7.4}. Statement (4) is a consequence of  Lemma \ref{5.3}. Statement (5) follows from (4) and Theorem \ref{6.7}(2). Statement (6) is obtained  from (2) and \ref{6.17}(1). Finally, statements (7) and (8) follow from (4)  together with  Theorem \ref{6.17}(2) and \ref{6.17}(1), respectively.
\end{proof}

\begin{lemma}\label{7.12}In $r_p(D(H_{\mathcal{D}}))$, we have
\begin{enumerate}
\item[(1)] $[V(l,\chi^{l-1})]=\sum_{i=0}^{[\frac{l-1}{2}]}(-1)^i\binom{l-1-i}{i}
x^i_{\phi\chi^2}y^{l-1-2i}$ for all $1\<l\<n$;
\item
[(2)]$[P(l,\chi^{l-1})]=\sum_{i=0}^{[\frac{n-l}{2}]}(-1)^i\frac{n-l}{n-l-i}
\binom{n-l-i}{i}x^{i+l-n}_{\phi\chi^2}y^{n-l-2i}[V(n,\chi^{n-1})]$ for all $1\<l<n$.
\end{enumerate}
\end{lemma}
\begin{proof}
(1) is similar to Lemma \ref{7.5}(3). We prove (2) by induction on $n-l$. If $l=n-1$, then  by Proposition \ref{7.11} (5), we have $[P(n-1,\chi^{n-2})]=x^{-1}_{\phi\chi^2}y[V(n,\chi^{n-1})]$ as desired. If $l=n-2$, then by Proposition \ref{7.11}(7),
$$\begin{array}{rl}
[P(n-2,\chi^{n-3})]=&x^{-1}_{\phi\chi^2}(y[P(n-1,\chi^{n-2})]-2[V(n,\chi^{n-1})])\\
=&(y^2x^{-1}_{\phi\chi^2}-2)x^{-1}_{\phi\chi^2}[V(n,\chi^{n-1})]\\
=&(y^2x^{-2}_{\phi\chi^2}-2x^{-1}_{\phi\chi})[V(n,\chi^{n-1})] \\
\end{array}$$
which matches the desired formula.

Now let $1\<l\<n-3$. Then by  Proposition \ref{7.11}(4),(8), together with the induction hypotheses, we obtain
$$\begin{array}{rl}
&[P(l,\chi^{l-1})]\\
=&x^{-1}_{\phi\chi^2}y[P(l+1,\chi^{l})]-x^{-1}_{\phi\chi^2}[P(l+2,\chi^{l+1})]\\
=&x^{-1}_{\phi\chi^2}y\sum_{i=0}^{[\frac{n-l-1}{2}]}(-1)^i\frac{n-l-1}{n-l-1-i}
\binom{n-l-1-i}{i}x^{i+l+1-n}_{\phi\chi^2}y^{n-l-1-2i}[V(n,\chi^{n-1})]\\
&-x^{-1}_{\phi\chi^2}\sum_{i=0}^{[\frac{n-l-2}{2}]}(-1)^i\frac{n-l-2}{n-l-2-i}
\binom{n-l-2-i}{i}x^{i+l+2-n}_{\phi\chi^2}y^{n-l-2-2i}[V(n,\chi^{n-1})]\\
=&\sum_{i=0}^{[\frac{n-l-1}{2}]}(-1)^i\frac{n-l-1}{n-l-1-i}
\binom{n-l-1-i}{i}x^{i+l-n}_{\phi\chi^2}y^{n-l-2i}[V(n,\chi^{n-1})]\\
&+\sum_{i=1}^{[\frac{n-l}{2}]}(-1)^i\frac{n-l-2}{n-l-1-i}
\binom{n-l-1-i}{i-1}x^{i+l-n}_{\phi\chi^2}y^{n-l-2i}[V(n,\chi^{n-1})].
\end{array}$$
If $n-l$ is odd, then
$$\begin{array}{rl}
[P(l,\chi^{l-1})]=&x^{l-n}_{\phi\chi^2}y^{n-l}+\sum_{i=1}^{[\frac{n-l}{2}]}(-1)^i
(\frac{n-l-1}{n-l-1-i}\binom{n-l-1-i}{i}\\
&+\frac{n-l-2}{n-l-1-i}
\binom{n-l-1-i}{i-1})x^{l+i-n}_{\phi\chi^2}y^{n-l-2i}\\
=&\sum_{i=0}^{[\frac{n-l}{2}]}(-1)^i\frac{n-l}{n-l-i}
\binom{n-l-i}{i}x^{i+l-n}_{\phi\chi^2}y^{n-l-2i}.
\end{array}$$

If $n-l$ is even, then
$$\begin{array}{rl}
[P(l,\chi^{l-1})]
=&x^{l-n}_{\phi\chi^2}y^{n-l}+\sum_{i=1}^{[\frac{n-l}{2}]-1}(-1)^i
(\frac{n-l-1}{n-l-1-i}\binom{n-l-1-i}{i}\\
&+\frac{n-l-2}{n-l-1-i}
\binom{n-l-1-i}{i-1})x^{l+i-n}_{\phi\chi^2}y^{n-l-2i}+
2(-1)^{\frac{n-l}{2}}x^{\frac{l-n}{2}}_{\phi\chi^2}\\
=&\sum_{i=0}^{[\frac{n-l}{2}]}(-1)^i\frac{n-l}{n-l-i}
\binom{n-l-i}{i}x^{i+l-n}_{\phi\chi^2}y^{n-l-2i}.
\end{array}$$
\end{proof}

\begin{corollary}\label{7.13}
The commutative ring $r_p(D(H_{\mathcal{D}}))$ is generated by $\{x_{\l},y, z_{\l'}|\l\in K, \l'\in I'_n\}$.
\end{corollary}
\begin{proof}
The proof  follows from Proposition \ref{7.11} and Lemma \ref{7.12}.
\end{proof}

\begin{proposition}\label{7.14}
In $r_p(D(H_{\mathcal{D}}))$, we have
$$\begin{array}{rl}
&(\sum_{i=0}^{[\frac{n}{2}]}(-1)^i\frac{n}{n-i}\binom{n-i}{i}x^{i-n}_{\phi\chi^2}
y^{n-2i}-\sum_{j=0}^{1}x^{-1}_{\phi^{jn}\chi^n})\\
&(\sum_{i=0}^{[\frac{n-1}{2}]}(-1)^i\binom{n-1-i}{i}x^i_{\phi\chi^2}y^{n-2i-1})=0\\
\end{array}$$
\end{proposition}
\begin{proof}
By Lemma \ref{7.12}(2), we obtain
 $$\begin{array}{rl}
 x^{-1}_{\phi\chi^2}yP(1,\varepsilon)
 =\sum_{i=0}^{[\frac{n-1}{2}]}(-1)^i\frac{n-1}{n-1-i}
\binom{n-1-i}{i}x^{i-n}_{\phi\chi^2}y^{n-2i}[V(n,\chi^{n-1})].
\end{array}$$

On the other hand, by Proposition \ref{7.11}(6) and Lemma \ref{7.12}(2), we have
$$\begin{array}{rl}
 &x^{-1}_{\phi\chi^2}yP(1,\varepsilon)\\
 =&x^{-1}_{\phi\chi^2}[P(2,\chi)]+(x^{-1}_{\chi^{n}}+x^{-1}_{\phi^n\chi^n})[V(n,\chi^{n-1})]\\
 =&\sum_{i=0}^{[\frac{n-2}{2}]}(-1)^i\frac{n-2}{n-2-i}
\binom{n-2-i}{i}x^{i+1-n}_{\phi\chi^2}y^{n-2-2i}[V(n,\chi^{n-1})]\\
&+(x^{-1}_{\chi^n}+x^{-1}_{\phi^n\chi^n})[V(n,\chi^{n-1})].\\
\end{array}$$
Therefore, we  have
$$\begin{array}{rl}
&\sum_{i=0}^{[\frac{n-1}{2}]}(-1)^i\frac{n-1}{n-1-i}
\binom{n-1-i}{i}x^{i-n}_{\phi\chi^2}y^{n-2i}
[V(n,\chi^{n-1})]-\\
&(\sum_{i=1}^{[\frac{n}{2}]}(-1)^{i-1}\frac{n-2}{n-1-i}
\binom{n-1-i}{i-1}x^{i-n}_{\phi\chi^2}y^{n-2i}+x^{-1}_{\chi^n}+x^{-1}_{\phi^n\chi^n})[V(n,\chi^{n-1})]=0.
\end{array}$$
Consequently, by Lemma \ref{7.12}(1), we obtain
$$\begin{array}{rl}
&(\sum_{i=0}^{[\frac{n}{2}]}(-1)^i\frac{n}{n-i}\binom{n-i}{i}x^{i-n}_{\phi\chi^2}
y^{n-2i}-\sum_{j=0}^{1}x^{-1}_{\phi^{jn}\chi^n})\\
&(\sum_{i=0}^{[\frac{n-1}{2}]}(-1)^i\binom{n-1-i}{i}x^i_{\phi\chi^2}y^{n-2i-1})=0.\\
\end{array}$$
\end{proof}
For any $1\<l<n$, let $f(x_{\phi\chi^2},y)=\sum_{i=0}^{[\frac{n-1}{2}]}(-1)^i\binom{n-1-i}{i}x^i_{\phi\chi^2}y^{n-1-2i}$ and
$$\begin{array}{rl}
&g_l(x_{\phi\chi^2},y)\\
=&(\sum_{i=0}^{[\frac{n-l}{2}]}(-1)^i\frac{n-l}{n-l-i}
\binom{n-l-i}{i}x^{i+l-n}_{\phi\chi^2}y^{n-l-2i})
(\sum_{i=0}^{[\frac{n-1}{2}]}(-1)^i
\binom{n-1-i}{i}x^{i}_{\phi\chi^2}y^{n-1-2i})\end{array}$$

\begin{proposition}\label{7.15}
Let $\l,\mu\in K$ and $\l',\l''\in I'_n$.
\begin{enumerate}
\item[(1)] $x_{\l}x_{\mu}=x_{\l\mu}.$
\item[(2)] $x_{\l}z_{\l'}=z_{\l\l'}.$
\item[(3)] $yz_{\l'}=z_{\chi\l'}+z_{\phi\chi\l'}.$
\item[(4)]  If $\l'\l''\in I'_n$, then $z_{\l'}z_{\l''}=\sum_{j=0}^{n-1}z_{\phi^j\l'\l''}.$
\item[(5)]  If $\l'\l''\in I_l$ for some $0\<l\<n-1$, then
    \end{enumerate}

    \begin{equation*}
\begin{array}{ll}
\vspace{0.2 cm}
z_{\l'}z_{\l''}=\left\{
\begin{array}{ll}
M_{l,\l',\l''}+f(x_{\phi\chi^2},y)x_{\phi^{\frac{n+l}{2}}\l'\l''\chi^{1-n}},& {\rm if }\ l \ {\rm and}\  n \ {\rm are\  odd},\\
M_{l,\l',\l''}, & {\rm if}\ l\ {\rm is \ odd}\ {\rm and} \ n \ {\rm is \ even},\\
M_{l,\l',\l''}+f(x_{\phi\chi^2},y)x_{\phi^{\frac{l}{2}}\l'\l''\chi^{1-n}},& {\rm if} \ l \ {\rm is \ even} \ {\rm and}  \ n \ {\rm is \ odd},\\
M_{l,\l',\l''}+\sum_{i=0}^1f(x_{\phi\chi^2},y)x_{\phi^{\frac{in+l}{2}}\l'\l''\chi^{1-n}}, & {\rm if}\ l \ {\rm and} \ n \ {\rm are \ even},\\
\end{array}\right.\\
\end{array}
\end{equation*}
where
$$\begin{array}{rl}
M_{l,\l',\l''}=& \sum_{c(l+1)\<s\<l}g_{n+l-2s}(x_{\phi\chi^2},y)x_{\phi^s\l'\l''\chi^{2s+1-l-n}}\\
  & +\sum_{c(n+l+1)\<s\<n-1}g_{2n+l-2s}(x_{\phi\chi^2},y)x_{\phi^s\l'\l''\chi^{1+2s-l-2n}}.\\
\end{array}$$
\end{proposition}
\begin{proof}
Statements (1) and (2) follow directly  from Lemma \ref{5.3}. Statement (3) follows from Theorem \ref{6.7}(1). Statement (4) is obtained  from Theorem \ref{6.16}(1). Statement (5) follows from Theorem \ref{6.16}(2), Proposition \ref{7.11} and Lemma \ref{7.12}.
\end{proof}

\begin{corollary}\label{7.16} The set
$$\{x_{\l}y^i,z_{\l'}|0\<i\<2n-2,\l\in K,\l'\in I'_n\}$$ forms a $\mathbb{Z}$-basis of $r_p(D(H_{\mathcal{D}}))$.
\end{corollary}
\begin{proof}
By Lemma \ref{7.14}, we have
$$\begin{array}{rl}
y^{2n-1}=&-\sum_{i=1}^{[\frac{n-1}{2}]}(-1)^i\binom{n-1-i}{i}x^i_{\phi\chi^2}y^{2n-2i-1}\\
&-\sum_{i=1}^{[\frac{n}{2}]}(-1)^i\frac{n}{n-i}\binom{n-i}{i}x^i_{\phi\chi^2}y^{2n-2i-1}+(x_{\phi^n\chi^n}+x_{\chi^n})y^{n-1}\\
&-(\sum_{i=1}^{[\frac{n}{2}]}(-1)^i\frac{n}{n-i}\binom{n-i}{i}x^i_{\phi\chi^2}y^{n-2i}-x_{\phi^n\chi^n}-x_{\chi^n})\\
&(\sum_{i=1}^{[\frac{n-1}{2}]}(-1)^i\binom{n-1-i}{i}x^i_{\phi\chi^2}y^{n-2i-1}).
\end{array}$$
It then follows from Proposition \ref{7.15} that $r_p(D(H_{\mathcal{D}}))$ is generated, as $\mathbb{Z}$-module,  by
$$\{x_{\l}y^i,z_{\l'}|0\<i\<2n-2,\l\in K,\l'\in I'_n\}.$$
Since  the rank of $\mathbb{Z}$-module $r_p(D(H_{\mathcal{D}}))$ is equal to  the cardinality of this set,  the set  $\{x_{\l}y^i,z_{\l'}|0\<i\<2n-2,\l\in K,\l'\in I'_n\}$ forms a  $\mathbb{Z}$-basis of $r_p(D(H_{\mathcal{D}}))$.
\end{proof}

Now we  describe the projective class ring $r_p(D(H_{\mathcal{D}}))$ in the two case of even $n$ and odd $n$, respectively. Let
$$Y=\{x_{\l},y, z_{\l'}|\l\in K, \l'\in I'_n\}$$
 and let $\mathbb{Z}[Y]$ denote the corresponding polynomial ring.

We first consider the case when $n$ is even.  Define a subset $J_0\subset \mathbb{Z}[Y]$ by
$$J_0=\left\{\left.\begin{array}{l}
x_{\l}x_{\mu}-x_{\l\mu};\  x_{\l}z_{\l'}-z_{\l\l'};\\
yz_{\l'}-z_{\chi\l'}-z_{\phi\chi\l'};\\
z_{\l_1}z_{\l_2}-\sum_{j=0}^{n-1}z_{\phi^j\l_1\l_2};\\
z_{\l_3}z_{\l_4}-M_{l,\l_3,\l_4};\\
z_{\l_5}z_{\l_6}-M_{k,\l_5,\l_6}-\\
\sum_{i=0}^1f(x_{\phi\chi^2},y)x_{\phi^{\frac{in+k}{2}}\l_5\l_6\chi^{1-n}};\\
(\sum_{i=0}^{\frac{n}{2}}(-1)^i\frac{n}{n-i}\binom{n-i}{i}x^{i-n}_{\phi\chi^2}y^{n-2i}\\
-\sum_{j=0}^1x^{-1}_{\phi^{jn}\chi^n})(\sum_{i=0}^{\frac{n-2}{2}}
(-1)^i\\
\binom{n-1-i}{i}x^i_{\phi\chi^2}y^{n-2i-1}).\\
 \end{array}
\right|\begin{array}{l}
\l,\mu\in K;\\
\l',\l_1,\l_2\in I'_n
{\rm with} \ \l_1\l_2\in I'_n;\\
\l_3,\l_4\in I'_n
{\rm with} \ \l_3\l_4\in I_l\\
{\rm for \ some \ odd} \ l , 0\<l\<n-1;\\
\l_5,\l_6\in I'_n
{\rm with} \ \l_5\l_6\in I_k\\
{\rm for \ some \ even} \ k, 0\<k\<n-1;\\
\end{array}\right\},$$
where $M_{l,\l_3,\l_4}$ and $M_{k,\l_5,\l_6}$ are given as in Proposition \ref{7.15}.

\begin{theorem}\label{7.17}  Suppose that $n$ is even. Then the projective class ring $r_p(D(H_{\mathcal{D}}))$ is isomorphic to the quotient ring  $\mathbb{Z}[Y]/( J_0)$ of $\mathbb{Z}[Y]$, where the ideal $(J_0)$ is generated by $J_0$.
\end{theorem}
\begin{proof}
The proof is similar to  that of Theorem \ref{7.10}.
\end{proof}

Now suppose that $n$ is odd. Define  a subset $J_1\subset \mathbb{Z}[Y]$ by

$$J_1=\left\{\left.\begin{array}{l}
x_{\l}x_{\mu}-x_{\l\mu};\  x_{\l}z_{\l'}-z_{\l\l'};\\
yz_{\l'}-z_{\chi\l'}-z_{\phi\chi\l'};\\
z_{\l_1}z_{\l_2}-\sum_{j=0}^{n-1}z_{\phi^j\l_1\l_2};\\
z_{\l_3}z_{\l_4}-M_{l,\l_3,\l_4}-\\
f(x_{\phi\chi^2},y)x_{\phi^{\frac{n+l}{2}}\l_3\l_4\chi^{1-n}};\\
z_{\l_5}z_{\l_6}-M_{k,\l_5,\l_6}-\\
f(x_{\phi\chi^2},y)x_{\phi^{\frac{k}{2}}\l_5\l_6\chi^{1-n}};\\
(\sum_{i=0}^{\frac{n-1}{2}}(-1)^i\frac{n}{n-i}\binom{n-i}{i}x^{i-n}_{\phi\chi^2}y^{n-2i}\\
-\sum_{j=0}^1x^{-1}_{\phi^{jn}\chi^n})(\sum_{i=0}^{\frac{n-1}{2}}
(-1)^i\\
\binom{n-1-i}{i}x^i_{\phi\chi^2}y^{n-2i-1}).\\
 \end{array}
\right|\begin{array}{l}
\l,\mu\in K;\\
\l',\l_1,\l_2\in I'_n
{\rm with} \ \l_1\l_2\in I'_n;\\
\l_3,\l_4\in I'_n
{\rm with} \ \l_3\l_4\in I_l\\
{\rm for \ some \ odd} \ l , 0\<l\<n-1;\\
\l_5,\l_6\in I'_n
{\rm with} \ \l_5\l_6\in I_k\\
{\rm for \ some \ even} \ k, 0\<k\<n-1;\\
\end{array}\right\},$$
where $M_{l,\l_3,\l_4}$ and $M_{k,\l_5,\l_6}$ are given as in proposition 6.16.

\begin{theorem}\label{nodd+1} Suppose that $n$ is odd. Then the projective class ring $r_p(D(H_{\mathcal{D}}))$ is isomorphic to the quotient ring  $\mathbb{Z}[Y]/( J_1)$ of $\mathbb{Z}[Y]$, where the ideal $(J_1)$ is generated by $J_1$.
\end{theorem}
\begin{proof}
The proof is analogous  to that of  Theorem \ref{7.10}.
\end{proof}

\begin{remark}
For any Hopf algebra $H$, the forgetful functor $F: {^H_H\mathcal{YD}} \to {_H\mathcal{M}}$, from the category of Yetter-Drinfeld $H$-modules to the category of left $H$-modules, is exact and monoidal.  Consequently, $F$ induces a ring homomorphism $\varphi$ from the Green ring $r({^H_H\mathcal{YD}})$ to the Green ring $r(H)$. Note that $r({^H_H\mathcal{YD}})$ is commonly denoted as $r(D(H))$ since $^H_H\mathcal{YD} \cong {_{D(H)}\mathcal{M}}$ holds when $H$ is finite-dimensional.

Similarly, the forget functor induces a  ring  homomorphism $\overline{\varphi}$ from the Grothendieck ring $G_0(D(H))$ to the Grothendieck ring $G_0(H)$.

It is clear that $\varphi$ (respectively, $\overline{\varphi}$) is surjective  if and only if every indecomposable (respectively, simple) $H$-module is the underlying module of some Yetter-Drinfeld $H$-module structure.

While this surjectivity holds for specific examples, such as the Taft algebras and the Radford algebras, this conclusion fails for a general pointed rank one Hopf algebra. Consequently, the study of the images of $\varphi$ and $\overline{\varphi}$ is of particular interest and warrants further investigation.
\end{remark}

\begin{example}
We consider the generalized Taft Algebra $H_{4,2}$.
Let $G=\langle g \rangle$ be a cyclic group of order $4$ generated by $g$, and let $q\in \Bbbk$ be a primitive $4$-th root of unit. Define a character $\gamma$ of $G$ by $\gamma(g)=q$. Then the character group is $\Gamma:={\rm Hom}(G,\Bbbk^{\times})=\{\gamma^i|i\in \mathbb{Z}_4\}$, and the dual group of $G\times \Gamma$ is $\Lambda=\widehat{G \times \Gamma}=\{\gamma^i\widehat{g}^j|i,j\in \mathbb{Z}_4\}$, where $\widehat{g}\in \widehat{\Gamma}$ is defined by $\widehat{g}(\gamma)=\gamma(g)=q$. Let $\chi=\gamma^2$ and $\rho=q^2$. Then $\mathcal{D}=(G,\chi,g,\a)$ is a group datum and $\chi(g)=\rho$, where $\a\in \Bbbk$. By a straightforward computation, one gets that
$$K=I_1=\{\gamma^{2j}\widehat{g}^j|j\in \mathbb{Z}_4\}, \ I_2=\{\gamma^{2(j-1)}\widehat{g}^j,\gamma\widehat{g}^j,\gamma^3\widehat{g}^j|j\in \mathbb{Z}_4\}.$$

Now let $\mathcal{D}=(G,\chi,g,0)$. This is a group datum of nilpotent type. The corresponding Hopf algebra $H_\mathcal{D}$ is the $8$-dimensional generalized Taft algebra $H_{4,2}$, which is generated, as an algebra, by $g$ and $x$ subject to the relations:
$$g^4=1, \quad x^2=0, \quad xg=\rho gx.$$The coalgebra structure of $H_{4,2}$ is given by $\Delta(g)=g\otimes g$ and $\Delta(x)=x\otimes g+1\otimes x$. By \cite{Liandzhang2013}, there are $4$ non-isomorphic simple modules over $H_{4,2}$:
$$\{M(1,\epsilon), M(1,\gamma), M(1,\gamma^2), M(1,\gamma^3)\},$$ and they are all one dimensional. On the other hand, all one-dimensional simple $D(H_{4,2})$-modules are $\{V(1,\lambda)|\lambda \in I_1\}$. It is obvious that for any $i$, $\gamma\widehat{g}^i \notin I_1$ and $\gamma^3\widehat{g}^i \notin I_1$. Consequently, the simple modules $M(1,\gamma)$ and $M(1,\gamma^3)$ do not possess Yetter-Drinfeld module structures.  It follows that the ring homomorphism $\overline{\varphi}: G_0(D(H_{4,2}))\rightarrow G_0(H_{4,2})$ is not surjective.

\end{example}

\vskip 12pt
\centerline{ACKNOWLEDGMENTS}

This work is supported by NNSF of China (Nos. 12201545, 12071412).\\

\end{document}